\documentclass[mathpazo,online]{jml}

\renewcommand\thisvolume{1}
\renewcommand\thisnumber{2}
\renewcommand\thisyear {2022}
\renewcommand\thismonth{June}

\renewcommand\originalreceiveddate{xx/xx/2022}

\renewcommand\acceptdate{xx/xx/2022}

\renewcommand\pagestart{1}
\renewcommand\pageend{xx}

\renewcommand\articledoi{doi.org/10.4208/jml.xxx}

\renewcommand\handler{xxx}

\usepackage{lipsum}

\usepackage[utf8]{inputenc}
\usepackage{mathrsfs,amsthm,xcolor,verbatim,bbm,amsmath,amsfonts,
	amssymb,nicefrac,enumitem}
\usepackage{hyperref,bm,mathtools,xparse,etoolbox}
\usepackage[capitalise,sort]{cleveref} 

\usepackage{algorithm}
\usepackage{algpseudocode}
\usepackage{multirow}
\usepackage{textcomp}
\usepackage{sidecap}
\pdfoutput=1 

\newcommand{\bx}{\boldsymbol{x}}
\newcommand{\bu}{\boldsymbol{u}}
\newcommand{\bdf}{\boldsymbol{f}}
\newcommand{\bg}{\boldsymbol{g}}

\begin{document}

	\title[RFM for Solving PDEs]{Bridging Traditional and Machine Learning-based Algorithms for Solving PDEs: The Random Feature Method}

	\author[1,2]{
		Jingrun Chen
		\thanks{
			{\tt jingrunchen@ustc.edu.cn}.
		}
	}
	\author[1]{
		Xurong Chi
		\thanks{
			{\tt cxr123@mail.ustc.edu.cn}.
		}
	}
	\author[3,4]{
		Weinan E
		\thanks{
			{\tt weinan@math.pku.edu.cn}.
		}
	}
	\author[1,4]{
		Zhouwang Yang
		\thanks{
			{\tt yangzw@ustc.edu.cn}.
		}
	}
	\affil[1]{School of Mathematical Sciences, University of Science and Technology of China}
	\affil[2]{Suzhou Institute for Advanced Research, University of Science and Technology of China}
	\affil[3]{AI for Science Institute, Beijing and Center for Machine Learning Research and School of Mathematical Sciences, Peking University}
	\affil[4]{School of Data Science, University of Science and Technology of China}
	
	\begin{abstract}
		One of the oldest and most studied subject in scientific computing is algorithms for solving partial differential equations (PDEs). A long list of numerical methods have been proposed and  successfully used for various applications. In recent years, deep learning methods have shown their superiority for high-dimensional PDEs where traditional methods fail. However, for low dimensional problems,  it remains unclear whether these methods have a real advantage over traditional algorithms as a direct solver. In this work, we propose the random feature method (RFM) for solving PDEs, a natural bridge between traditional and machine learning-based algorithms. RFM is based on a combination of well-known ideas: 1. representation of the approximate solution using random feature functions; 2.  collocation method to take care of the PDE; 3. the penalty method to treat the boundary conditions, which allows us to treat the boundary condition and the PDE in the same footing. We find it crucial to add several additional components including multi-scale representation and rescaling the weights in the loss function. We demonstrate that the method exhibits spectral accuracy and can compete with  traditional solvers in terms of both accuracy and efficiency. In addition, we find that RFM is  particularly suited for complex problems with complex geometry, where both traditional and machine learning-based algorithms encounter difficulties.
	\end{abstract}
	
	%
	%
	%
	\keywordone{Partial differential equation,}
	\keywordtwo{Machine learning,}
	\keywordthree{Random feature method,}
	\keywordfour{Rescaling}
	
	\maketitle
	
	\section{Introduction}
	\label{sec1}
	
	One of the oldest and most studied subject in scientific computing is algorithms for solving partial differential equations (PDEs).
	Finite difference \cite{leveque2007finite}, finite element \cite{zienkiewicz2005finite}, spectral methods \cite{shen2011spectral} and a host of other methodologies have been proposed and studied, with great success.
	At the same time,  a variety of scientific softwares based on these methodologies have been developed and widely used by the academia
	as well as the industry. They have become standard resources in most, if not all, engineering applications.
	
	In recent years, as neural network models have had great success in a variety of artificial intelligence (AI) tasks, the idea of using these models 
	to solve PDEs has gained a lot of popularity \cite{EHJ2017, HJE2018, EY2018, SS2018, PINN, Bao}.
	Though back in the 90's, it was already proposed to use neural networks as test or trial functions in PDE solvers \cite{LeeKang1990}, 
	the recent proposals often have some non-trivial new twist. 
	The most notable success is to solve PDEs
	and control problems in high dimensions \cite{EHan2016, EHJ2017, HJE2018, SS2018}, a class of problems that 
	traditional algorithms are not able to handle. 
	Indeed  deep learning-based algorithms have now made it fairly routine to solve
	a large class of PDEs in hundreds of  or even higher dimensions \cite{review}, something  impossible to do just a few years ago.
	In another direction, neural networks can also be used to parametrize the solution operator of PDEs  \cite{Khoo, DeepOnet, Fouriernet}, which is also beyond the capability of traditional algorithms.
	
	
	Despite these great deal of efforts and the great deal of success, the situation with solving PDEs is not entirely satisfactory
	even for some of the traditional engineering problems.  Here is an incomplete list of some of the difficulties we still encounter.
	
	\begin{enumerate}
		\item Problems with complex geometry. A typical problem is Stokes flow in porous media \cite{allen2021mathematics}.
		In principle the finite element method (FEM) is ideally suited for problems with complex geometry. In practice, coming up with
		a suitable mesh is often a highly non-trivial task both in terms of the human efforts  and the actual computational cost required.
		Machine learning-based algorithms, while easy to code, have not proven to be reliable and competitive in practical situations against
		traditional algorithms.
		
		\item Kinetic equations. Although its dimensionality is much lower than the high dimensional ones mentioned above,
		kinetic equations such as the Boltzmann equation are traditionally regarded as high dimensional problems for which
		classical methods do encounter difficulties. Ideas based on sparse grids should help \cite{Griebel, ShenJie}, but at the moment the most
		popular approach for solving kinetic equations is still the direct simulation Monte Carlo algorithm (DSMC) \cite{stefanov2019basic}.
		One problem with DSMC is that the solutions produced contain too much noise.
		
		\item Multi-scale problems. Examples include
		problems involving chemical kinetics that typically span a large range of time scales; fully developed turbulent flows
		that contain a large range of spatial and temporal scales; and the modeling of composite materials; see \cite{weinan2011principles} for example.
	\end{enumerate}
	
	
	
	
	
	Our objective in this paper is to propose a methodology for solving general PDEs that shares the merits of both classical and
	machine learning-based algorithms.
	This new class of algorithms can be made spectrally accurate. At the same time, they are also mesh-free,
	making them easy to use even in settings with complex geometry.
	Our starting point is based on a combination of rather simple and well-known ideas: 
	We use random feature functions to represent the approximate solution, the collocation method to take care of the PDE as well as the boundary conditions in the least-squares sense, and a rescaling procedure to balance the contributions from the PDE and the boundary conditions in the loss function. In actual implementations, we take several inspirations from the machine learning literature.
	\begin{enumerate}
		\item Our preferred choices of the basis functions are random feature functions. 
		As a result, the method bears close similarity
		with the random feature model in machine learning.
		Though deterministic basis functions can also work, we have found that random feature functions 
		are generally more reliable and more efficient.
		If necessary, the actual choice of the basis functions can be tuned beforehand in a precomputing stage 
		in order to make them more adapted to the nature of the problem. 
		For this reason, we will 
		name  the class of algorithms proposed here  the``random feature method (RFM)''.
		
		The feature (or basis) functions adopted here are the ones used in neural networks (for details, see below).
		This means that we are using a special class of random feature models: models that come from two-layer neural networks
		with the inner parameters fixed. As a comparison, if we use radial basis functions as the feature functions, we are uncertain 
		whether the performance of the RFM would be equally good. We also find that the additional gain in training the inner parameters
		does not offset at all the increased complexity in the training.
		
		Another important component is a multi-scale representation of the solutions. We use a partition of unity (PoU) to piece 
		together different local representations as well as a global representation for the large scale components of the solution.
		This strategy has proven to be vital in practice in order to achieve good accuracy. 
		
		\item A rescaling procedure is needed to balance the contributions from the PDE and the boundary/initial conditions in the loss function,
		by tuning the weight parameters ($\{\lambda_{Ii}^{k}\}$ and $\{\lambda_{Bj}^{\ell}\}$ in \eqref{loss2}).
		Although the situation is similar to the one in training neural network models, the reduced complexity from using the random feature model
		instead of the neural network model seems to make the task of parameter tuning much simpler. In fact, this rescaling procedure
		is quite simple and fully automatic.
		
	\end{enumerate}
	
	A closely related method is the local extreme learning machine (locELM) proposed in \cite{dong2021local}, using the combination of extreme learning machines \cite{huang2006extreme} and domain decomposition. 
	The idea of domain decomposition can be viewed as a special choice of PoU.
	In that case, some smoothness conditions have  to be enforced explicitly across the boundary of the different domains. 
	For a series of one-dimensional and two-dimensional problems with simple geometries and explicit solutions, this method shows spectral accuracy.  
	However, it does not seem to work well for more practical problems such as the linear elasticity problem, even with simple geometry. 
	
	This paper is organized as follows. In Section \ref{sec2}, we present the RFM: the construction of approximate solution, the loss function, and the optimization procedure. This is followed by results of two sets of  numerical experiments: 
	one with explicit solutions and the other with complex geometries. 
	We use the former to demonstrate that RFM has spectral accuracy, and  the latter to demonstrate its feasibility for solving complex problems. 
	These results are shown in Section \ref{sec3}. In Section \ref{sec4}, we present some discussions.
	
	\section{The random feature method}
	\label{sec2}
	
	Consider the following problem 
	\begin{equation}
	\begin{cases}
	\mathcal{L}\bu(\boldsymbol{x}) =  \bdf(\boldsymbol{x})  & \boldsymbol{x}\in\Omega,\\
	\mathcal{B}\bu(\boldsymbol{x}) =  \bg(\boldsymbol{x})  & \boldsymbol{x}\in \partial\Omega,
	\end{cases}
	\label{pde}
	\end{equation}
	where $\boldsymbol{x}=(x_{1},\cdots,x_{d})^{T}$, and $\Omega$ is bounded and connected domain in $\mathbb{R}^d$. Examples include the elliptic problem, the linear elasticity problem, and the Stokes flow problem.
	
	Roughly speaking, much like the random feature model in machine learning, RFM relies on three key components: 1. The loss function is built on the least-squares (strong) formulation of the PDEs on collocation points; 2. The approximate solution is constructed using a set of random feature functions; 3. The training is very much like neural network training, with the additional step of rescaling the penalty parameters to balance the contributions from different terms. In what follows, we will discuss each component in some detail.
	
	\subsection{Loss function}
	\label{sec2-1}
\begin{figure}[htbp]
		\centering
		\includegraphics[width=0.4\textwidth]{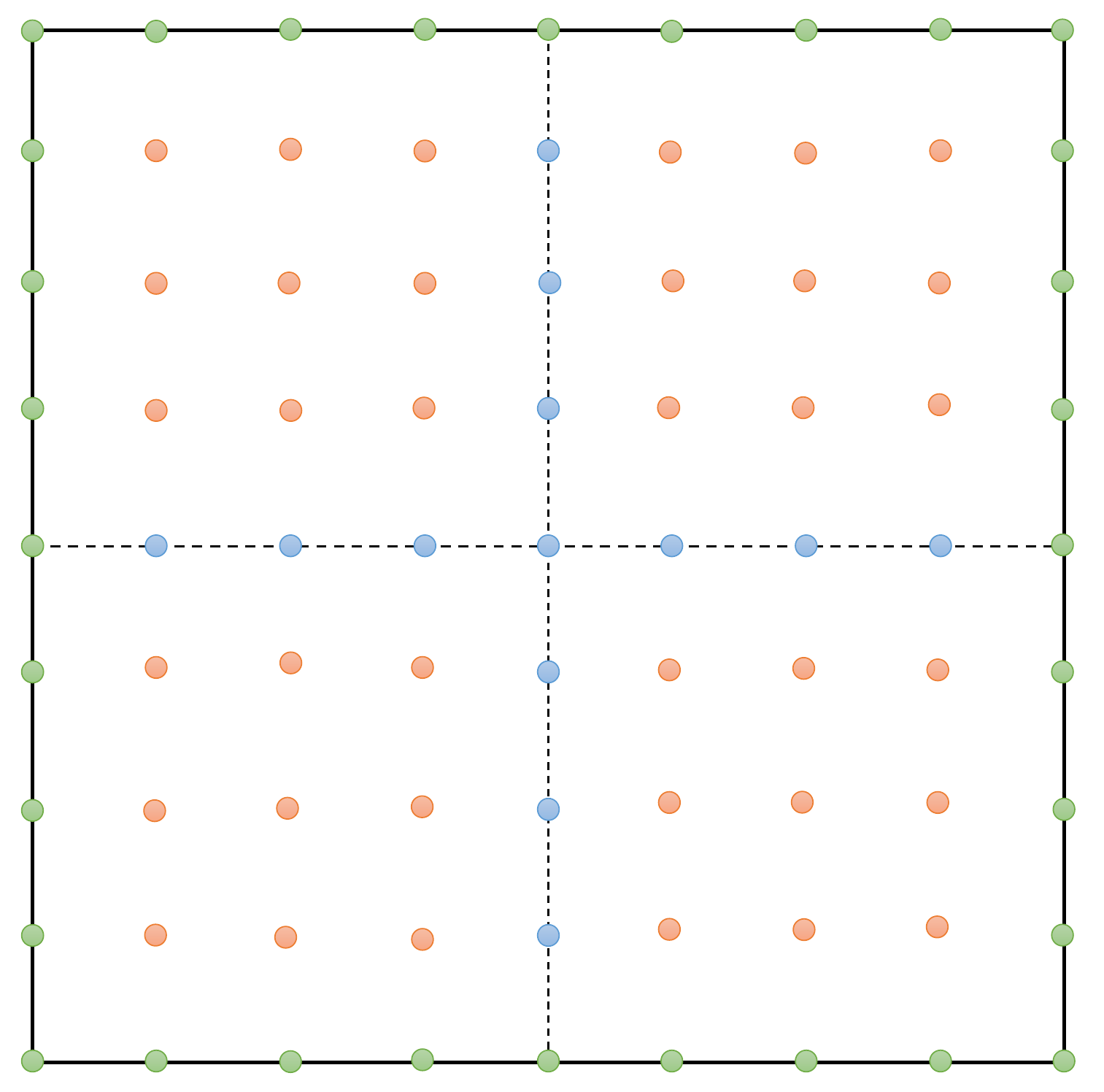}
		\caption{Collocation points for a square domain: $C_{I}$, interior points in orange and blue; $C_{B}$, boundary points in green.}
		\label{collocation_points}
	\end{figure}

	There are three standard approaches for solving \eqref{pde}: the weak form,   the strong form and the variational form in cases when 
	\eqref{pde} is the Euler-Lagrange equation of some variational problem.
	Each of these approaches gives rise to some particular choices of loss function.
	For neural network-based algorithms, examples of these different loss functions can be found in \cite{EY2018,PINN,Bao}.
	In this paper, we will focus on the strong form at collocation points to construct the loss function. Corresponding to \eqref{pde}, we have two sets of collocation points: $C_{I}$,  the set of interior points in $\Omega$ and $C_{B}$, the set of boundary points on $\partial\Omega$. 
	Let $C= C_{I}\cup C_{B}$ be the set of all collocation points. 
	At each collocation point, we will enforce either the PDE or the boundary condition.
	Let $K_I$ and $K_B$ be the number of conditions at each interior point and boundary point, respectively.
	The total number of conditions is $N = K_I \# C_{I} + K_B \# C_{B}$.
	See Figure \ref{collocation_points} for an illustration. Detailed selection algorithm for $C$ will be specified later.
	
	A simple choice of the loss function for \eqref{pde} is as follows:
	\begin{equation}
	Loss = \sum_{\boldsymbol{x}_{i} \in C_{I}}\sum_{k=1}^{K_I}\lambda_{Ii}^{k}\|\mathcal{L}^{k}\bu(\boldsymbol{x}_{i})-\bdf^{k}(\boldsymbol{x}_{i})\|_{l^{2}}^{2}+ \sum_{\boldsymbol{x}_{j} \in C_{B}}\sum_{\ell=1}^{K_B}\lambda_{Bj}^{\ell}\|\mathcal{B}^{\ell}\bu(\boldsymbol{x}_{j})-\bg^{\ell}(\boldsymbol{x}_{j})\|_{l^{2}}^{2}.
	\label{loss2}
	\end{equation}
	Here $\{\lambda_{Ii}^k\}$ and $\{\lambda_{Bj}^{\ell}\}$ are the penalty parameters.  
	In this form, we allow different choices of the penalty parameters at different collocation points.
	By treating the boundary conditions and the PDE in the same footing, we do not need to impose boundary conditions for the
	feature function. This gives us much needed flexibility for treating problems with complex geometry.
	
	\subsection{Random feature functions}
	\label{sec2-2}
	Following the random feature model in machine learning, we construct the approximate solution $u_M$ of $u$ by a linear combination of $M$ network basis functions $\{\phi_m\}$ over $\Omega$ as follows
	\begin{equation}\label{eqn:rfm}
	u_M(\boldsymbol{x}) = \sum_{i=1}^M u_m \phi_m(\boldsymbol{x}).
	\end{equation}
	For  vectorial solutions, we approximate each component of the solution using \eqref{eqn:rfm}, i.e. 
	$$\bu_M(\bx) = (\sum_{i=1}^M u_m^1\phi_{m}^1(\boldsymbol{x}),\cdots, \sum_{m=1}^M u_m^{K_I}\phi_{m}^{K_I}(\boldsymbol{x}))^T.$$
	
	Generally speaking the basis functions will be chosen as the ones that occur naturally in neural networks, for example:
	$$
	\phi_m(\boldsymbol{x}) = \sigma(\boldsymbol{k}_m \cdot\boldsymbol{x} + b_m),
	$$
	where $\sigma$ is some scalar nonlinear function, $\boldsymbol{k}_m, b_m $ are some random but fixed parameters.
	For solving PDE problems, activation functions such as $\tanh$, $\sin$, and $\cos$ can all be used.
	
	

	In practice, additional ideas are needed to achieve good performance.
	
	\subsubsection{Partition of unity and local random feature models}\label{sec2-2-1}
	Random feature functions are globally defined, while the solution of the PDE typically has local variations, possibly at small scales.
	To accommodate this, we construct many local solutions, each of which corresponds to a random feature model, and piece them
	together using partition of unity (PoU).
	
	To construct the PoU, 
	we start with a set of points $\{\boldsymbol{x}_n\}_{n=1}^{M_p}\subset\Omega$, each of which serves as the center for a component
	in the partition. 
	For each $n$, construct the normalized coordinate:
	\begin{equation}\label{eqn:normalization}
	\tilde{\boldsymbol{x}}=\frac{1}{\boldsymbol{r}_{n}}(\boldsymbol{x}-\boldsymbol{x}_{n}), \quad n=1,\cdots, M_p,
	\end{equation}
	where $\boldsymbol{r}_n = (r_{n1}, r_{n2}, \cdots, r_{nd})$ and $\{ \boldsymbol{r}_n\}$ is preselected. This linear transformation 
	maps $[x_{n1}-r_{n1},x_{n1}+r_{n1}]\times \cdots \times [x_{nd}-r_{nd},x_{nd}+r_{nd}]$ onto $[-1,1]^{d}$. 
	
	Next, we construct $J_n$ random feature functions by
	\begin{equation}\label{eqn:basis0}
	\phi_{nj}(\boldsymbol{x}) = \sigma(\boldsymbol{k}_{nj} \cdot \tilde{\boldsymbol{x}} + b_{nj}), \quad j=1, \cdots, J_n,
	\end{equation}
	where the feature vectors $\{(\boldsymbol{k}_{nj},  b_{nj})\}$  often chosen randomly. 
	A common choice is the uniform distribution $\boldsymbol{k}_{nj} \sim \mathbb{U}([-R_{nj},R_{nj}]^{d})$ and $b_{nj} \sim \mathbb{U}([-R_{nj},R_{nj}])$, though different distributions can be used. In this way the locally space-dependent information is incorporated into $M= \sum_{n=1}^{M_p} J_n$ random feature functions.
	
\begin{figure}[htbp]
		\centering
		\includegraphics[width=0.6\textwidth]{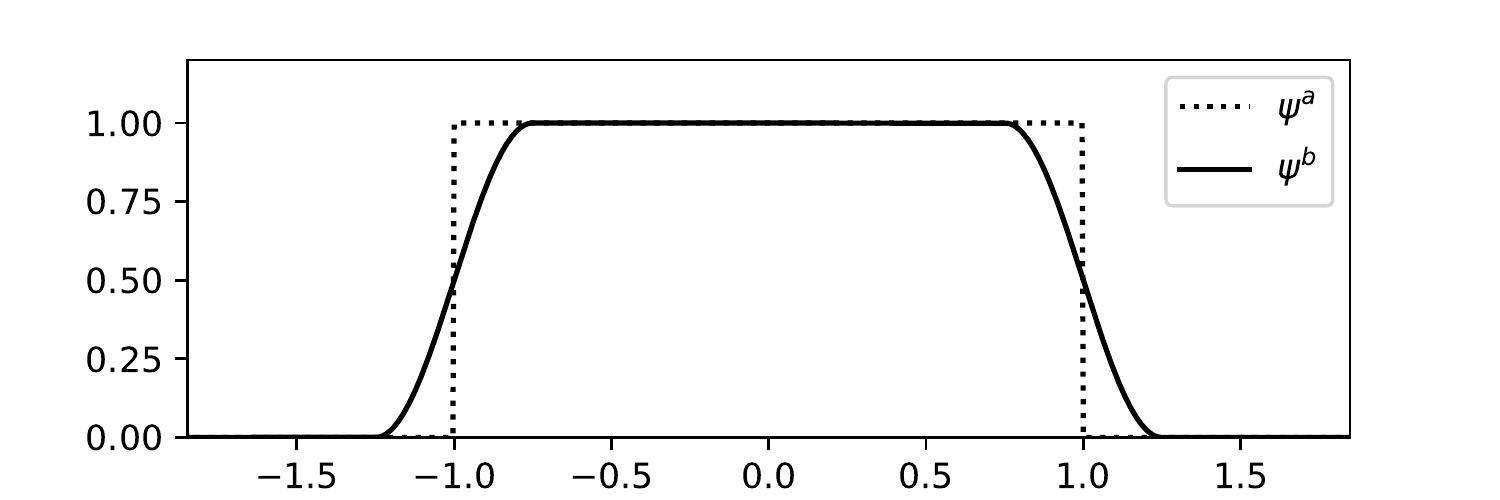}
		\caption{Visualization of $\psi^{a}(x)$ in \eqref{psi1} and $\psi^{b}(x)$ in \eqref{psi2}.}
		\label{PoU}
	\end{figure}

	We now discuss the construction of the PoU.
	When $d=1$,  let 
	\begin{equation}
	\psi_{n}^{a}(x)=\mathbb{I}_{-1 \leq \tilde{x} < 1},
	\label{psi1}
	\end{equation}
	and
	\begin{equation}\label{psi2}
	\psi_{n}^{b}(x) = 	\left\{\begin{aligned}
	&\frac{1+\sin(2\pi \tilde{x})}{2} && -\frac{5}{4}\leq \tilde{x} < -\frac{3}{4}, \\
	&1 && -\frac{3}{4}\leq \tilde{x} < \frac{3}{4}, \\
	&\frac{1-\sin(2\pi \tilde{x})}{2} && \frac{3}{4}\leq \tilde{x} < \frac{5}{4}, \\
	&0 && \text{otherwise}.
	\end{aligned}\right.
	\end{equation}
	See Figure \ref{PoU} for the visualization of $\psi^{a}$ and $\psi^{b}$.
	
	High-dimensional PoU can be constructed using the tensor product of one-dimensional PoU functions $\psi_n(\boldsymbol{x})=\prod \limits_{k=1}^{d}\psi_n(x_{k})$. 
	
	Putting together, the approximate solution $u_M$ in \eqref{eqn:rfm} is given by
	\begin{equation}
	u_M(\boldsymbol{x})=\sum_{n=1}^{M_p} \psi_n (\boldsymbol{x})   \sum_{j=1}^{J_n }u_{nj} \phi_{nj} (\boldsymbol{x}).
	\label{representation2}
	\end{equation}

	\subsubsection{Multi-scale basis}\label{sec2-2-2}
	In some situations, \eqref{representation2} alone is less efficient in capturing the large scale features in the solution.
	Therefore, on top of the PoU-based local basis functions, we can add another global component:
	\begin{equation}
	u_M(\boldsymbol{x})= u_g(\boldsymbol{x})+  \sum_{n=1}^{M_p} \psi_n (\boldsymbol{x})   \sum_{j=1}^{J_n }u_{nj} \phi_{nj} (\boldsymbol{x})
	\label{representation3}
	\end{equation}
	where $u_g$ is a global random feature function; see \eqref{eqn:rfm}. 
	
	\subsubsection{Adaptive basis}\label{sec2-2-3}
	The ideal choice of the distribution for the feature vectors is one that reflects the spectral distribution of the solution,
	which is not available to us beforehand.
	In some situations, we can obtain some incomplete information about the spectral distribution of the solution
	in the precomputing stage.
	For example, if the PDE has an inhomogeneous forcing term, we can perform a spectral analysis of the forcing term.
	The result can be used to guide us in the selection of the spectral distribution of the feature vectors.
	We have seen that this is particularly useful when $\sin$/$\cos$ is used as the activation function.
	
	Such a procedure can be seen as a compromise between the random feature method and the two-layer neural network model.
	In the two-layer neural network model, the inner parameters, i.e. the feature vectors, are part of the parameter set to be optimized.
	This in principle allows us to obtain an optimal choice of the feature vectors.
	The price we pay is that the optimization problem is much more complicated than simply fixing the feature vectors.
	In the random feature model, the feature vectors are fixed. If we choose a wrong set of feature vectors, the accuracy will deteriorate.
	With an adaptive procedure, one might be able to avoid this and at the same time still retain the simplicity of a linear model.
	

	\subsection{Optimization}\label{sec2-3}
	Recall the loss function:
	\begin{align}
	Loss  = \sum_{\boldsymbol{x}_{i} \in C_{I}}\sum_{k=1}^{K_I}\lambda_{Ii}^{k}\| \mathcal{L}^{k} \bu_M(\bx_i)-\bdf^{k}(\boldsymbol{x}_{i})\|_{l^{2}}^{2}  + \sum_{\boldsymbol{x}_{j} \in C_{B}}\sum_{\ell=1}^{K_B}\lambda_{Bj}^{\ell}\| \mathcal{B}^{\ell} \bu_M(\bx_j)-\bg^{\ell}(\boldsymbol{x}_{j})\|_{l^{2}}^{2}\label{loss3},
	\end{align}
	where 
	$$\bu_M(\bx) = (\sum_{n=1}^{M_p} \psi_n (\boldsymbol{x})   \sum_{j'=1}^{J_n }u^1_{nj'} \phi^1_{nj'} (\boldsymbol{x}),\cdots, \sum_{n=1}^{M_p} \psi_n (\boldsymbol{x})   \sum_{j'=1}^{J_n }u_{nj'}^{K_I} \phi^{K_I}_{nj'} (\boldsymbol{x}))^T.$$
	
	This optimization problem can be solved using standard algorithms for least-squares approximation.
	One important new twist is the tuning of the penalty parameters.
	One simple yet effective way is to rescale the penalty parameters so that each term in the loss function 
	is of the same order of magnitude.  This can be done using the largest term in the loss function as the reference.
	The detailed formula will be shown in the next section.
	This is particularly important in situations where the physical constants in the PDEs are of disparate size.
	
	\subsection{Collocation points}\label{sec2-4}
	There are many existing collocation point sampling methods for general geometric representations. If the boundary has an explicit parametric representation, the collocation points can be chosen as uniform grid points in the parameter space.
	In the case when the boundary has an implicit geometric representation, we can easily identify the interior points and define 
	an  energy function for finding a point on the boundary. In the implementations presented below,  we uniformly sample $Q$ points over a rectangle $R$ containing $\Omega$ and delete points in $R \cap \Omega^{c}$. 
	The construction of collocation points over a two-dimensional rectangular domain $\Omega$ is illustrated in Figure \ref{collocation_points}.
	
	
	We summarize the main steps of RFM in Algorithm \ref{alg:RFM}.
	\begin{algorithm}[htbp]
		\caption{The random feature method.}
		\label{alg:RFM}
		\begin{algorithmic}[1]
			\Statex \textbf{Input:}
			Number of basis functions $M$;
			number of collocation points $Q$;
			rule for generating collocation points;
			\Statex \textbf{Output:}
			The approximate solution $u_M$;
			\State Construct $M$ random feature functions $\{\phi_{m}\}$ and the PoU  $\{\psi_{n}\}$;
			\State Sample points $C= C_{I}\cup C_{B}$ according to some predetermined rule;
			\State Evaluate equations at $C_{I}$ and boundary conditions at $C_{B}$;
			\State Construct the loss function \eqref{loss3} ($M$ is not necessarily equal to $N$);
			\State Solve the optimization problem; 
			\State Return $u_M$;
		\end{algorithmic}
	\end{algorithm}
	
	
	\section{Numerical results}
	\label{sec3}
	
	We report results for two kinds of situations:  We use problems with explicit solutions to study how the performance of RFM depends on
	the different components in the algorithms such as the choice of the basis functions.
	We also present results for problems that do not have explicit solutions to demonstrate the power of RFM in complicated situations.
	
	The examples discussed below are all second-order PDEs.  We need  $C^1$ smoothness for the approximate solution. If we use $ \{\psi^{a}\}$, we impose this smoothness condition explicitly on the collocation points at the interfaces of the elements in the partition.
	If we use $ \{\psi^{b}\}$, no additional smoothness conditions need to be imposed.
	
	In what follows, unless indicated otherwise, we use the default setup where the collocation points are uniformly distributed on $\Omega$, 
	the weights $\{k_m\}$ and $\{b_m\}$ are assumed to follow the distribution $\mathbb{U}[-1,1]$, the activation function is chosen to be $\tanh$, 
	the PoU is $ \{\psi^{a}\}$. To achieve good accuracy,  we find that the distribution of $\{k_m\}$ and $\{b_m\}$ should be weakly problem-dependent when $\sin$ and/or $\cos$ is used as the activation function. $R_m=1$  works well for the examples discussed if $\tanh$ is used as the activation function. 
	We first select a set of points $\{\boldsymbol{x}_{n}\}_{n=1}^{M_p}$ and construct the PoU as follows. 
	For each $\boldsymbol{x}_{n}$, we construct $J_{n}$ random feature functions with radius $\boldsymbol{r}_{n}$. 
	Then we sample $Q$  equally spaced collocation points. The ones that are outside $\Omega$ are deleted.
	

	When evaluating the error, we take a refined grid with grid size being half of the size for the collocation points.
	The errors are evaluated on this refined grid.

	

	We start with some simple examples and use them to study how the performance of RFM, particularly the accuracy,
	depends on the details of algorithm.
	
	\subsection{Choice of random feature functions}
	\label{sec3-1}
	\begin{example}[Helmholtz equation] 
		Consider the one-dimensional Helmholtz equation with Dirichlet boundary condition over $\Omega= [0,8]$
		\begin{equation}
		\left\{ \begin{aligned}
		& \frac{\mathrm{d}^{2} u(x)}{\mathrm{d} x^{2}}-\lambda u(x) = f(x)  \quad x\in\ \Omega,\\
		&u(0)=c_{1}, \quad u(8)=c_{2}.
		\end{aligned} \right.
		\label{helmholtz}
		\end{equation}
		Once an explicit form of $u$ is given, $c_{1}$, $c_{2}$, and $f$ can be computed.
	\end{example}
	
	\begin{example}[Poisson equation] 
		Consider the Poisson equation with Dirichlet boundary condition over $\Omega= [0,1]\times[0,1]$
		\begin{equation}
		\left\{ \begin{aligned}
		& \Delta u(x, y)=f(x, y)  \quad (x,y)\in\Omega,\\
		& u(x, 0)=g_{1}(x), \quad u(x, 1)=g_{2}(x), \\
		& u(0, y)=h_{1}(y), \quad u(1, y)=h_{2}(y).
		\end{aligned} \right.
		\label{poisson}
		\end{equation}
		Again once an explicit form of $u$ is given, $g_{1}$, $g_{2}$, $h_{1}$, $h_{2}$, and $f$  can be computed.
	\end{example}
	
	Detailed numerical results for these problems can be found in Appendix \ref{sec::appendix::basis}.  Here we briefly summarize the main findings.
	
	We find that as long as the support of the distribution for the weights approximately covers the frequency domain of the true solution, RFM produces stable results using $\sin$/$\cos$ as the activation functions. In addition, we observe that in most cases, random sampling performs better than deterministic choices of the feature vectors.
	
	In addition, the introduction of multi-scale basis functions in \eqref{representation3} improves the  accuracy. 
	A Fourier analysis of the error confirms that  the use of multi-scale basis functions reduce the low-frequency error more effectively.
	
	We see also that with $\psi^{a}$,  the error  is more concentrated around the interface between the macro-elements in the partition.
	With $\psi^{b}$, the error  tends to be concentrated at the boundary.
	
	One point of interest is the comparison between RFM and PINN. We observe that the accuracy of PINN is around $1E-3$.
	Increasing the network size does not seem to improve the accuracy. In contrast,  we observe exponential rate of convergence for RFM in terms of the number of random feature functions.
	
	I  

	
	
	\subsection{Rescaling}
	\label{sec3-2}
	An important consideration is the balance between contributions from the PDE terms and the boundary conditions in the loss function. 
	This requires tuning the weights of different terms. 
	In this subsection, we consider the elasticity problem when $d=2$ and demonstrate how the rescaling strategy works.
	The idea is to rescale each term in the loss function to the same order of magnitude according to the largest term in the sum. Specifically, we choose the penalty parameters in \eqref{loss3} as follows:
	\begin{align}
	&\lambda_{Ii}^{k} = \frac{c}{\underset{1\leq n\leq M_p}{\max}\underset{1\leq j'\leq J_n}{\max}\underset{1\leq k'\leq K_I}{\max}|\mathcal{L}^{k} (\phi^{k'}_{nj'}(\boldsymbol{x}_{i})\psi_{n}(\boldsymbol{x}_{i}))|}\quad \boldsymbol{x}_{i} \in C_{I}, \; k = 1,\cdots, K_I, \label{penalty1}\\
	&\lambda_{Bj}^{\ell} = \frac{c}{\underset{1\leq n\leq M_p}{\max}\underset{1\leq j'\leq J_n}{\max}\underset{1\leq \ell'\leq K_I}{\max}|\mathcal{B}^{\ell} (\phi^{\ell'}_{nj'}(\boldsymbol{x}_{j})\psi_{n}(\boldsymbol{x}_{j}))|}\quad \boldsymbol{x}_{j} \in C_{B}, \; \ell = 1,\cdots, K_B, \label{penalty2}
	\end{align}
	where $c$ is a universal constant and we set $c=100$.
	
	We will see that this simple strategy significantly improves the accuracy, particularly in situations when the  physical constants in the PDE
	are of disparate size.
	
	The two-dimensional elasticity problem we consider here is of the following form
	\begin{equation}
	\left\{\begin{aligned}
	-\operatorname{div}(\boldsymbol{\sigma}(\boldsymbol{u}(\boldsymbol{x})))  &= \boldsymbol{B}(\boldsymbol{x}) \quad && \boldsymbol{x}\in\Omega, \\ 
	\boldsymbol{\sigma}(\boldsymbol{u}(\boldsymbol{x})) \cdot \boldsymbol{n}  &=  \boldsymbol{N}(\boldsymbol{x}) \quad && \boldsymbol{x}\in \Gamma_N,\\
	\boldsymbol{u}(\boldsymbol{x}) \cdot \boldsymbol{n}  &=  \boldsymbol{U}(\boldsymbol{x}) \quad && \boldsymbol{x} \in \Gamma_D,
	\end{aligned}\right.
	\label{elasticity}
	\end{equation}
	where $\boldsymbol{\sigma}: \mathbb{R}^{2} \rightarrow \mathbb{R}^{2}$ is the stress tensor induced by the displacement field $\boldsymbol{u}: \Omega \rightarrow \mathbb{R}^{2}$, $\boldsymbol{B}$ is the body force over $\Omega$, $\boldsymbol{N}$ is the surface force on $\Gamma_N$, $\boldsymbol{U}$ is the displacement on $\Gamma_D$, and $\partial\Omega=\Gamma_N\cup\Gamma_D$.
	
	Following \cite{nguyen2008meshless}, we consider the  Timoshenko beam problem with size $L \times D$,
	subject to a parabolic traction at the free end as shown in Figure \ref{fig3-2-1}. The exact solution and the experimental setup are presented in Appendix \ref{sec::appendix::complex::rescaling}.
	\begin{figure}[htbp]
		\centering
		\includegraphics[width=0.5\textwidth]{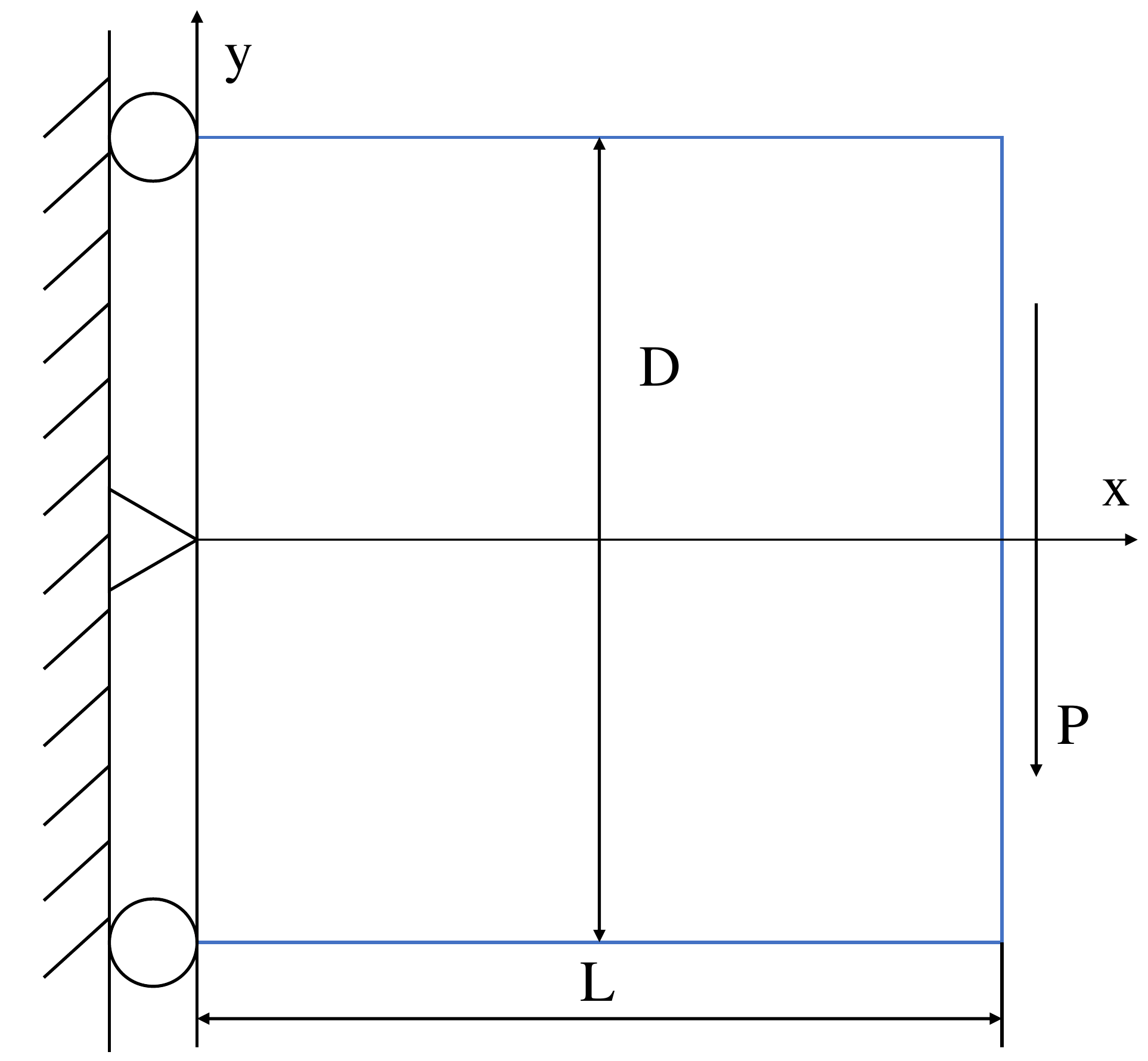}
		\caption{The Timoshenko beam problem.}
		\label{fig3-2-1}
	\end{figure}
	
	Let us first look at the accuracy of RFM and locELM. To this end, we presented in Table \ref{table3-2-1} the relative error.
	Since $\sigma_{y}$ is zero everywhere,  we omit this term. 
	From Table \ref{table3-2-1}, we observe that the relative error of locELM is around $1E-3$, while RFM still has  spectral accuracy. 
	We attribute the improved performance of RFM  to the rescaling strategy.
	\begin{table}[htbp]
		\caption{\label{table3-2-1} Comparison of RFM and locELM for the Timoshenko beam problem.} \centering
		\begin{small}
			\begin{tabular}{|c|cc|cccc|}
				\hline
				Method & $M$ & $N$ & $u$ error  & $v$ error & $\sigma_{x}$ error  & $\tau_{xy}$ error \\
				
				\hline
				\multirow{8}{*}{RFM}
				& \multirow{4}{*}{800} & 400   & 1.36E-2  & 3.43E-3  & 1.40E-2   & 1.63E-2  \\
				&                      & 1200  & 7.14E-6  & 7.98E-7  & 8.93E-6   & 7.45E-6  \\
				&                      & 4000  & 6.41E-11 & 4.34E-11 & 6.41E-11  & 6.58E-11 \\
				&                      & 14400 & 8.16E-12 & 1.01E-12 & 1.07E-11  & 1.03E-11 \\
				\cline{2-7}
				& \multirow{4}{*}{3200}& 1680  & 1.02E-2  & 1.42E-3  & 1.13E-2  & 7.65E-3  \\
				&                      & 4960  & 4.51E-6  & 7.89E-7  & 4.98E-6  & 4.36E-6  \\
				&                      & 16320 & 1.22E-11 & 7.23E-12 & 1.56E-11 & 1.40E-11 \\
				&                      & 58240 & 5.17E-13 & 1.49E-13 & 1.47E-12 & 1.99E-11 \\
				
				\hline
				\multirow{8}{*}{locELM}
				& \multirow{4}{*}{800} & 400   & 5.22E-3  & 4.90E-3  & 1.33E-2  & 2.39E-2 \\
				&                      & 1200  & 1.55E-4  & 5.25E-5  & 1.44E-4  & 1.02E-4 \\
				&                      & 4000  & 6.36E-4  & 3.47E-4  & 6.55E-4  & 7.26E-4 \\
				&                      & 14400 & 1.76E-3  & 1.64E-3  & 1.93E-3  & 2.57E-3 \\
				\cline{2-7}
				& \multirow{4}{*}{3200}& 1680  & 8.50E-2  & 4.04E-2  & 7.72E-2  & 4.19E-2 \\
				&                      & 4960  & 1.32E-5  & 6.19E-6  & 3.25E-5  & 4.22E-5 \\
				&                      & 16320 & 1.33E-3  & 1.12E-3  & 1.31E-3  & 1.04E-3 \\
				&                      & 58240 & 6.42E-4  & 1.91E-4  & 1.18E-3  & 1.38E-3 \\
				\hline
			\end{tabular}
		\end{small}
	\end{table}
	
\begin{figure}[htbp]
		\centering
		\includegraphics[width=0.8\textwidth]{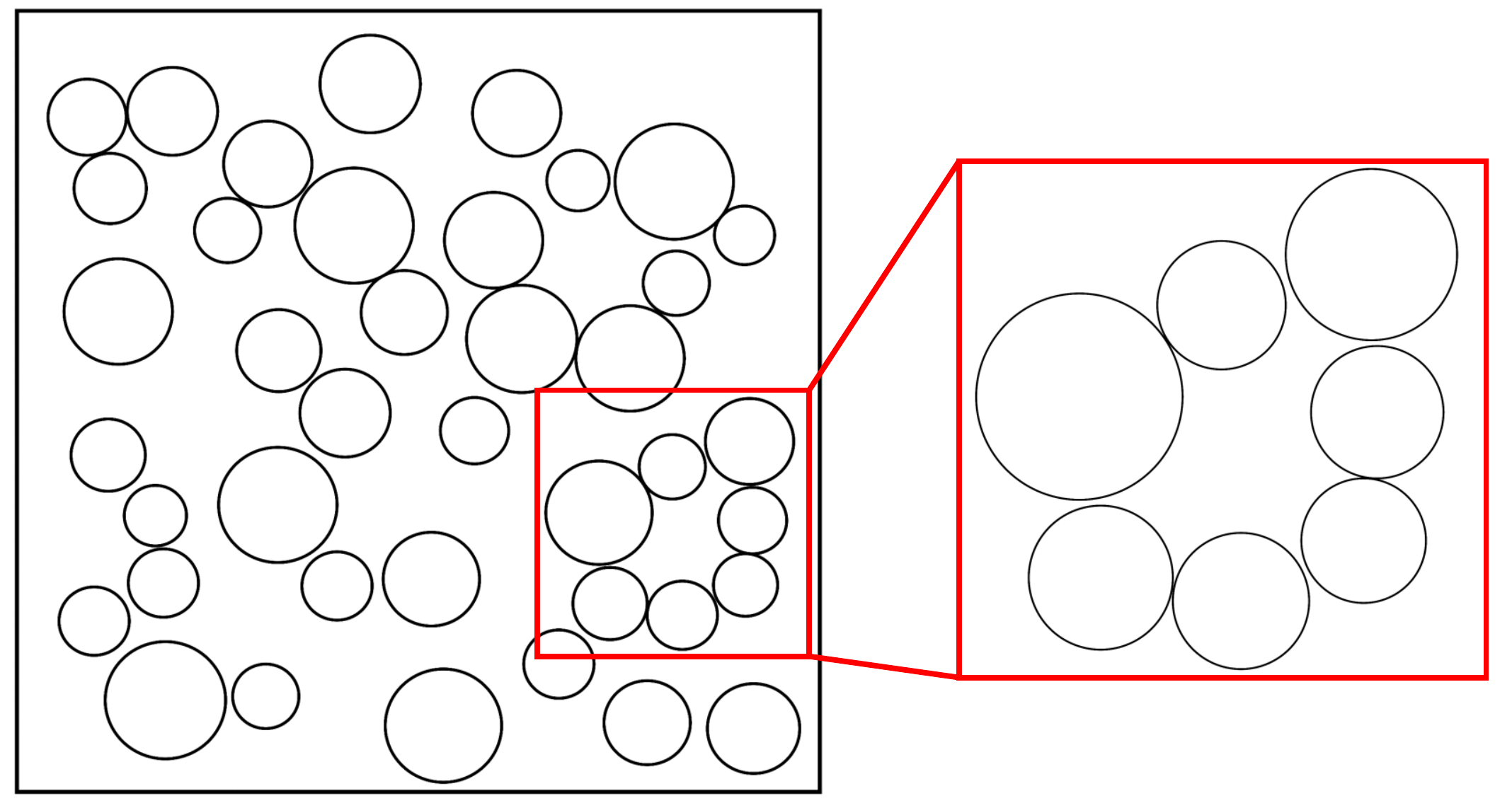}
		\caption{A two-dimensional complex domain.}
		\label{fig3-2-2}
	\end{figure}

	Next, we study	a two-dimensional elasticity problem with a complex geometry; see Figure \ref{fig3-2-2}. 
	Here $\Omega$ is defined as a square $(0,8) \times (0,8)$ with $40$ holes of radius between $0.3 $ and $0.6]$ inside.
	Note that there is a cluster of holes that are nearly touching, as shown in the inset. 
	More details can be found in Appendix \ref{sec::appendix::complex::rescaling}.
	
	In this example, the error produced by locELM is around $10^{-3}\sim10^{-2}$, while RFM still maintains spectral accuracy, as is shown in Table \ref{table3-2-2}, Appendix \ref{sec::appendix::complex::rescaling}.
	

	\subsection{Comparison with FEM}
	\label{sec3-3}
	In this subsection, we compare RFM with the classical adaptive FEM for  two elasticity problems. 
	
\begin{figure}[htbp]
		\centering
		\subfigure[$u$]{
			\includegraphics[width=0.45\textwidth]{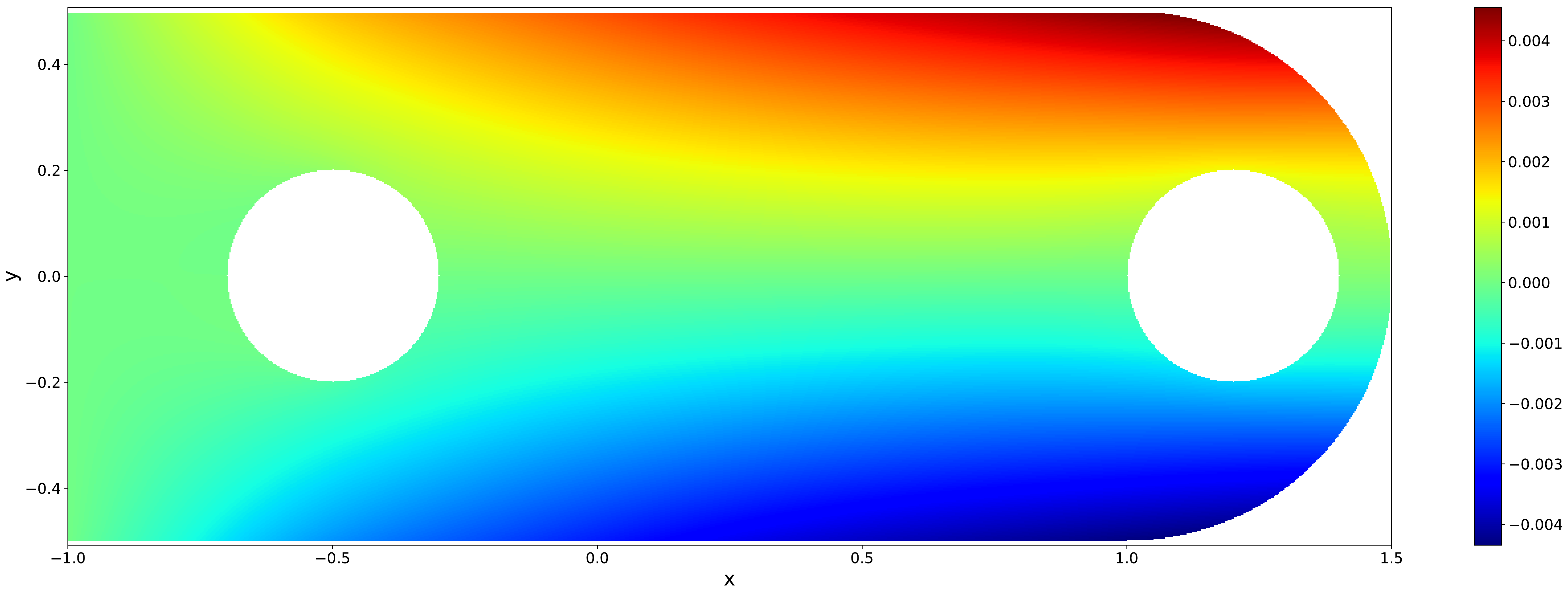}
		}
		\quad
		\subfigure[$v$]{
			\includegraphics[width=0.45\textwidth]{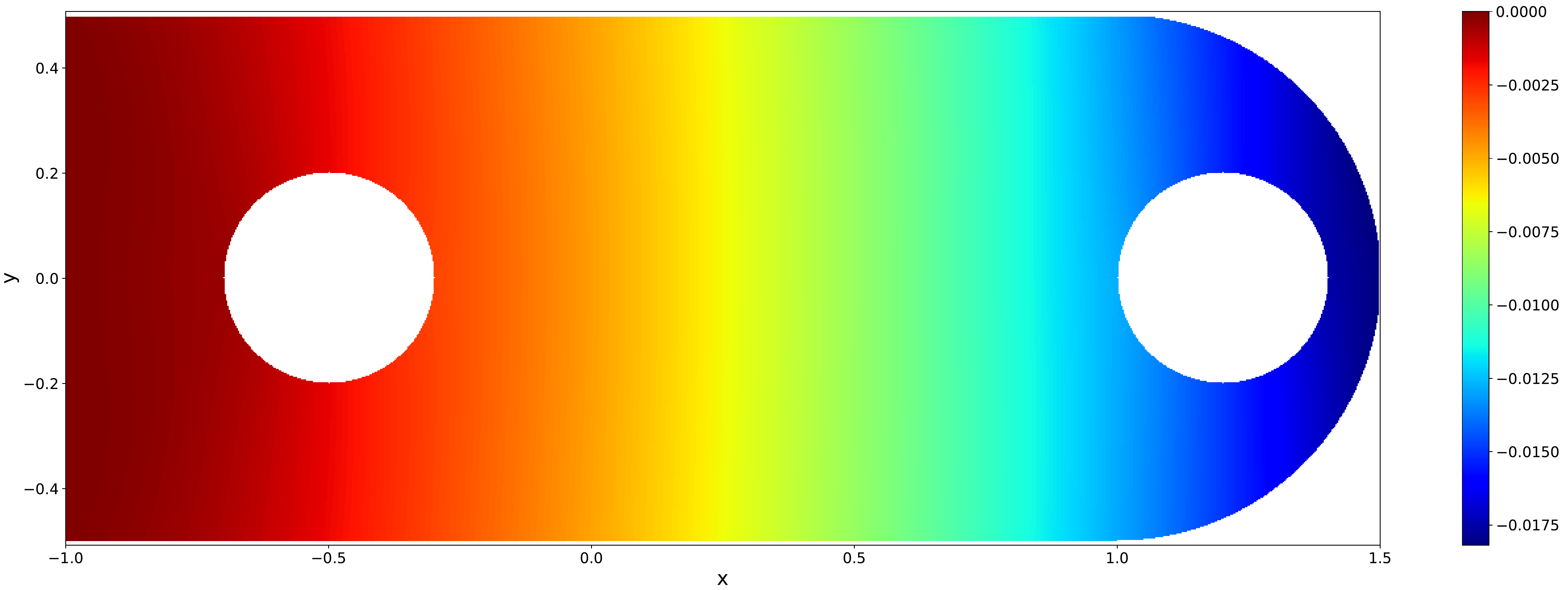}
		}
		\quad
		\subfigure[$\sigma_{x}$]{
			\includegraphics[width=0.28\textwidth]{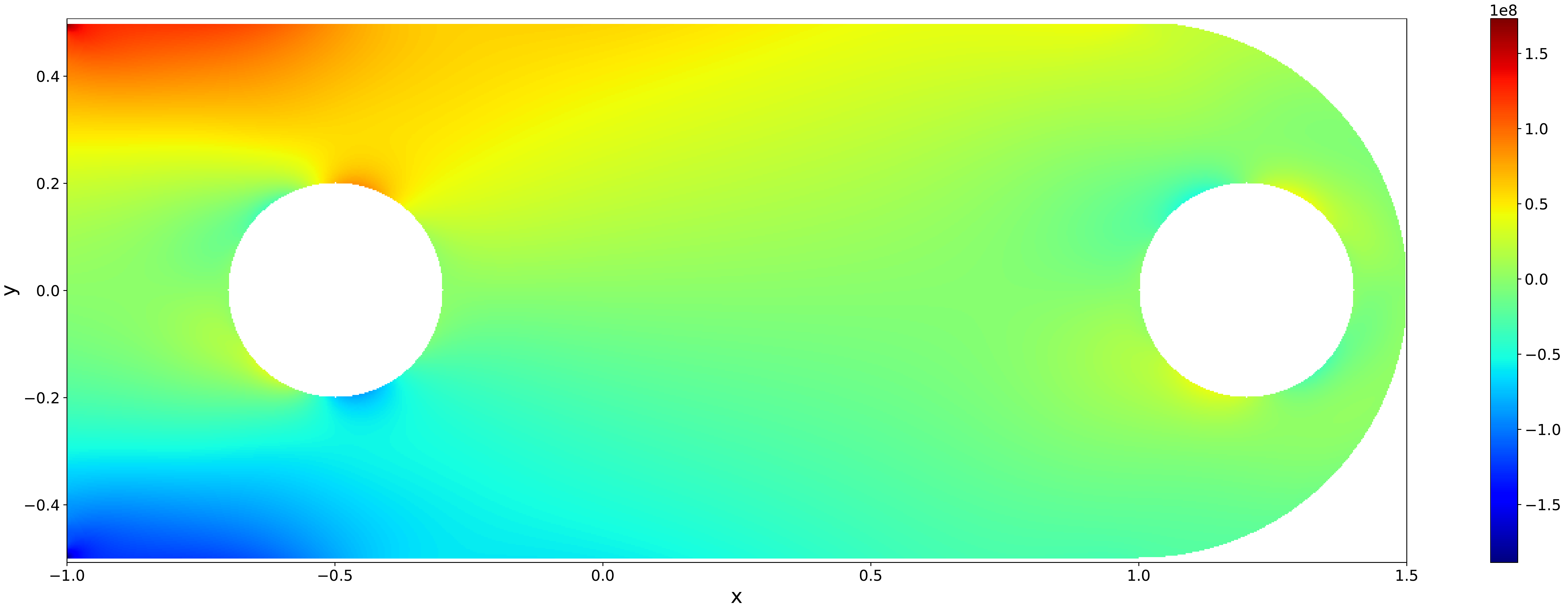}
		}
		\quad
		\subfigure[$\tau_{xy}$]{
			\includegraphics[width=0.28\textwidth]{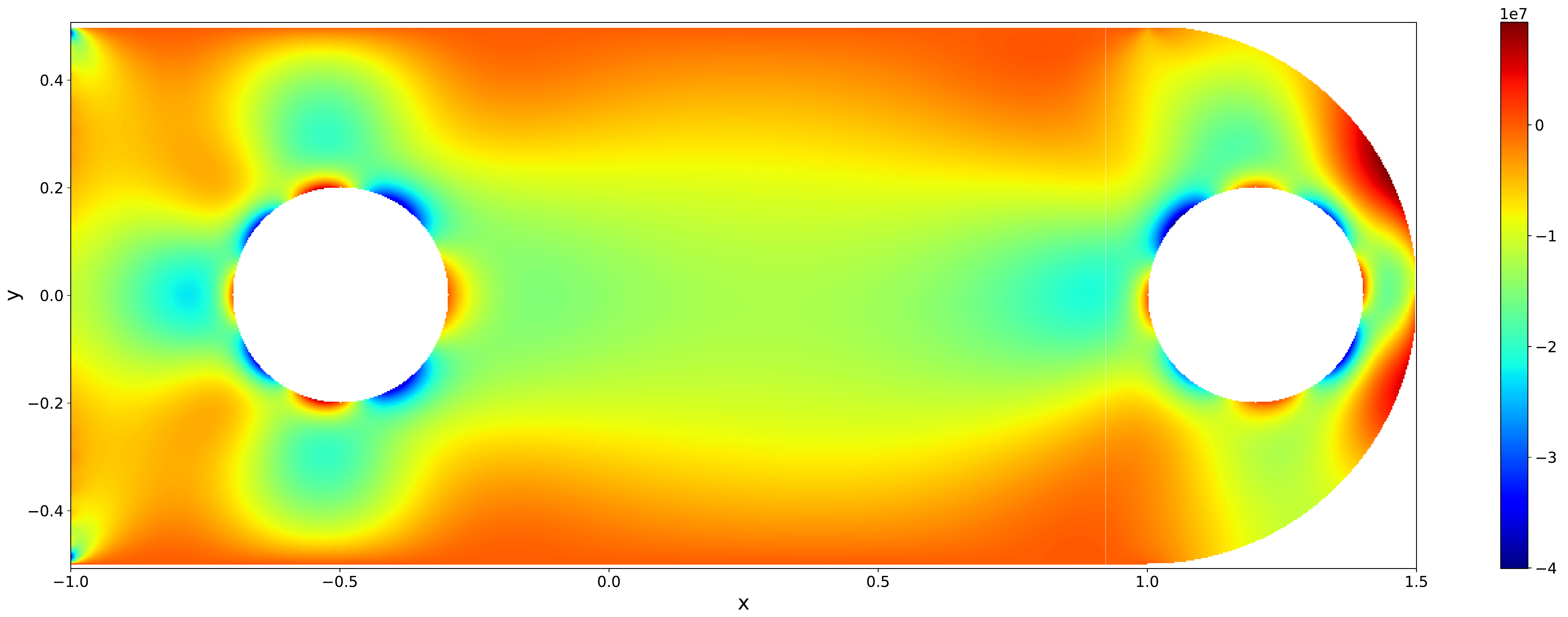}
		}
		\quad
		\subfigure[$\sigma_{y}$]{
			\includegraphics[width=0.28\textwidth]{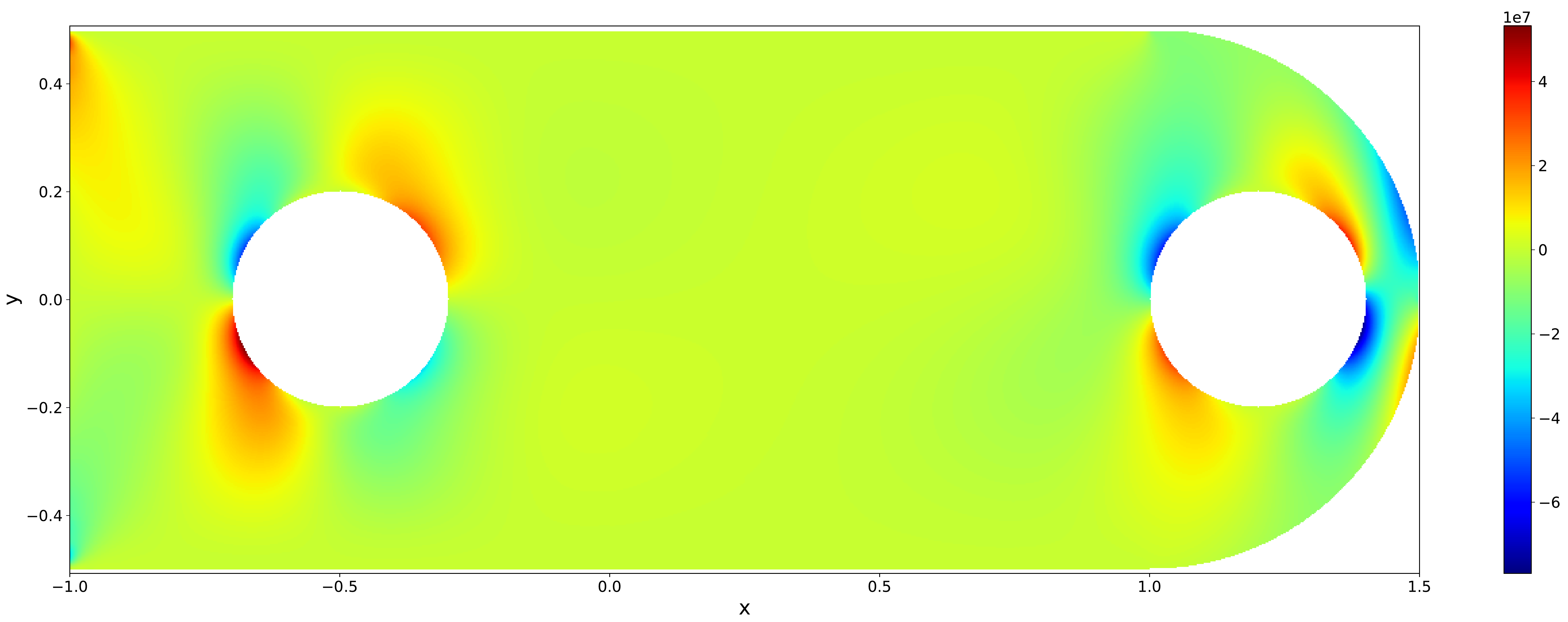}
		}
		\caption{Numerical solution by the random feature method for the  elasticity problem.}
		\label{fig3-3-1}
	\end{figure}

	The domain for the first example is given by a square $(-1,1) \times(-0.5,0.5)$ jointed by a semi-disk centered at $(1.0,0.0)$ with radius $0.5$,
	with two disks centered at $(1.2,0.0),(-0.5,0.0)$ with radius $0.2$ removed (see Figure \ref{fig3-3-1}).
	The material parameters are the same as those in Section \ref{sec3-2}. The left boundary $x=-0.5$ is fixed and a load $P=10^{7}$ Pa is applied on the upper half of the semicircle. Dirichlet boundary condition is applied on the left boundary $x=-0.5$ and Neumann boundary condition is applied on the other boundaries.
	
	More details of the experimental setup can be found in Appendix \ref{sec::appendix::complex::fem}. Figure \ref{fig3-3-1} visualizes the displacement fields $u$, $v$, and the stress fields $\sigma_x$, $\tau_{xy}$, $\sigma_y$.
	For this example, it is quite straightforward to obtain a FEM solution. 
	From Table \ref{table3-3-1} in Appendix \ref{sec::appendix::complex::fem}, we see that the difference between the RFM and FEM  solutions is about $1\%$.

	For the second example, we use  the same domain as in Figure \ref{fig3-2-2}, Section \ref{sec3-2} and the same materials constants. The lower boundary $y=0$ is fixed, and a load of $10^{5}\sin (x+y) e^{y}$ Pa along the positive $x$ direction is applied on both the left and right boundaries. 
	More details of the experimental setup can be found  in Appendix \ref{sec::appendix::complex::fem}.
	
	\begin{figure}[htbp]
		\centering
		\subfigure[$u$]{
			\includegraphics[width=0.45\textwidth]{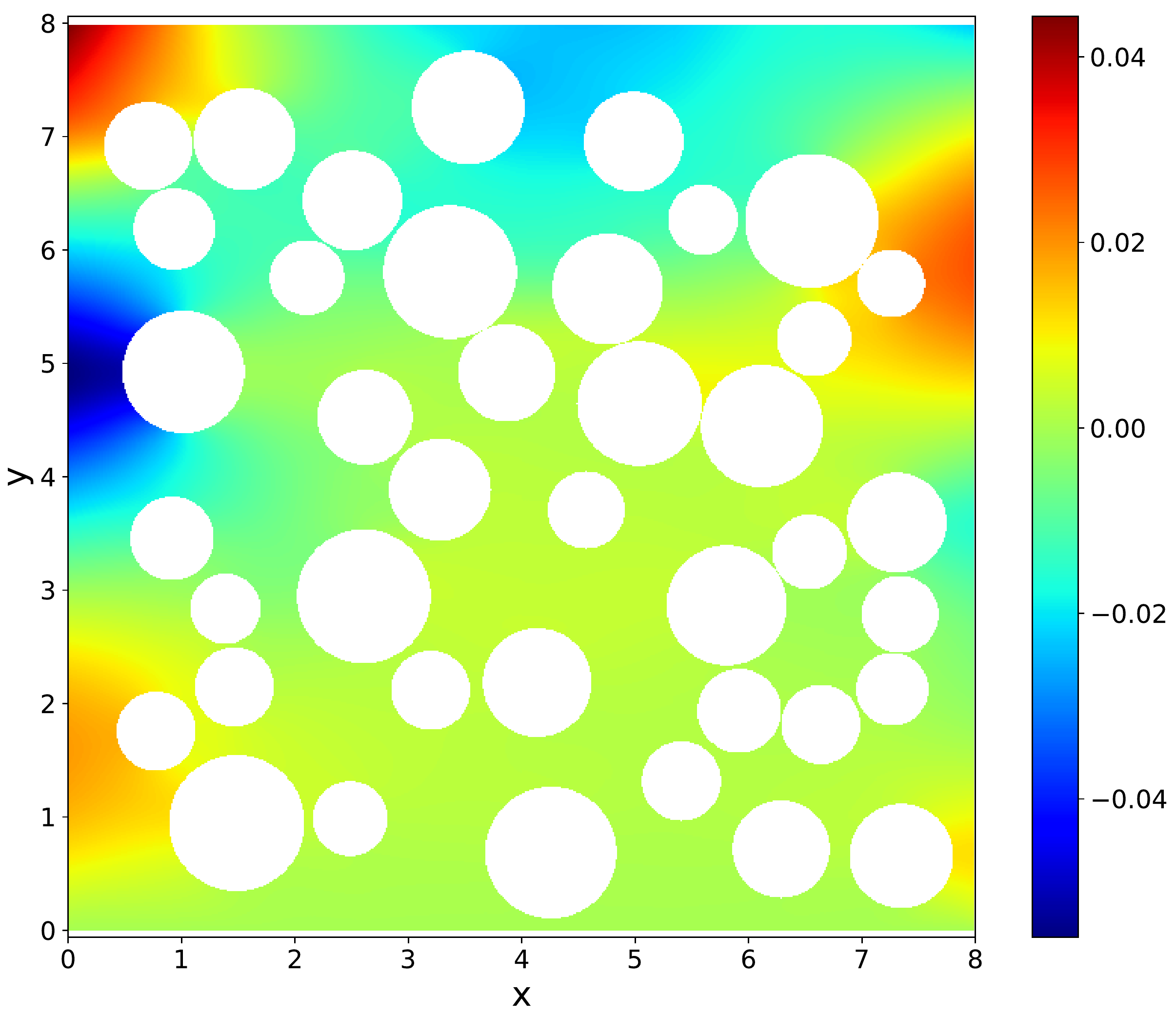}
		}
		\quad
		\subfigure[$v$]{
			\includegraphics[width=0.45\textwidth]{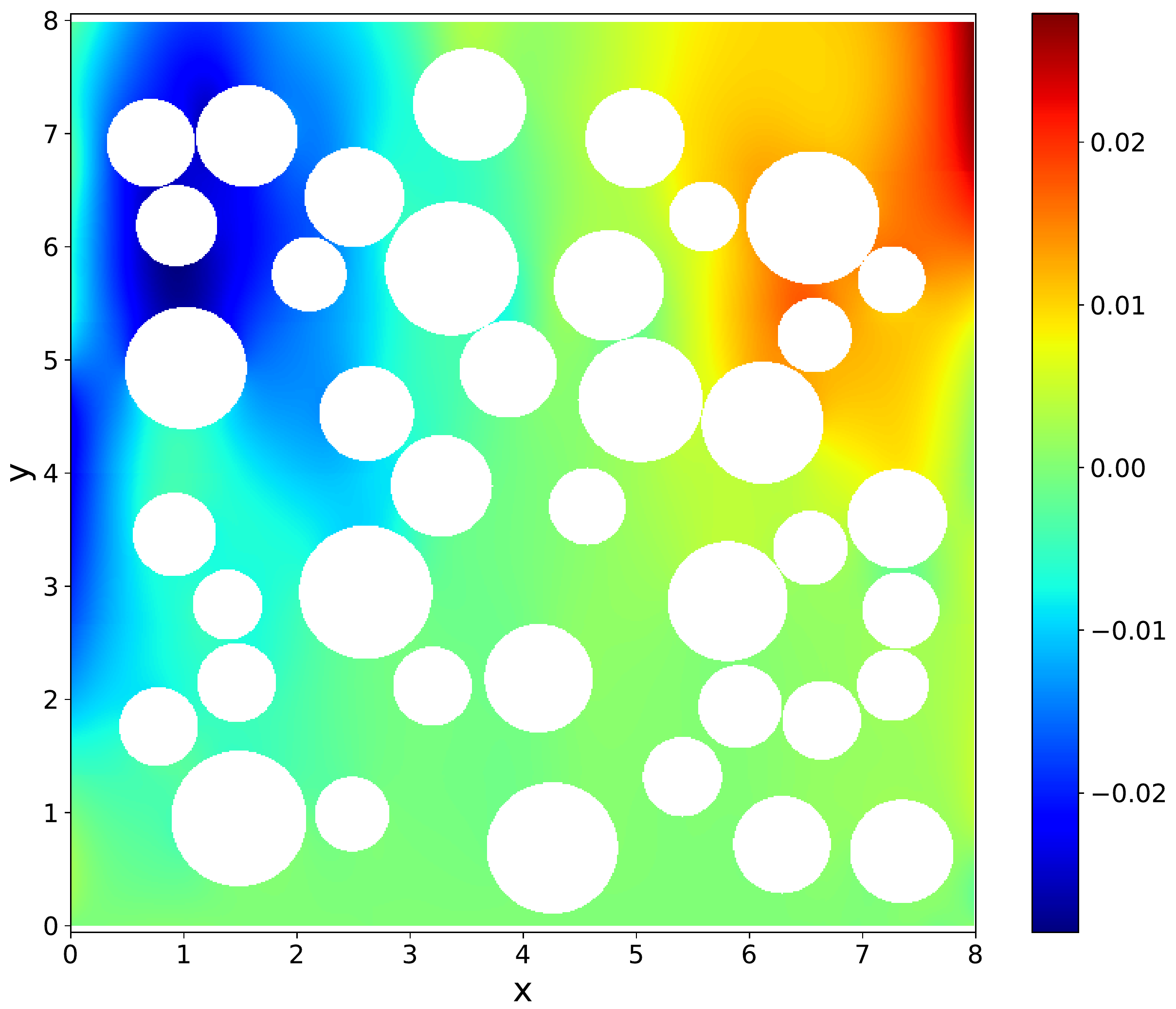}
		}
		\quad
		\subfigure[$\sigma_{x}$]{
			\includegraphics[width=0.45\textwidth]{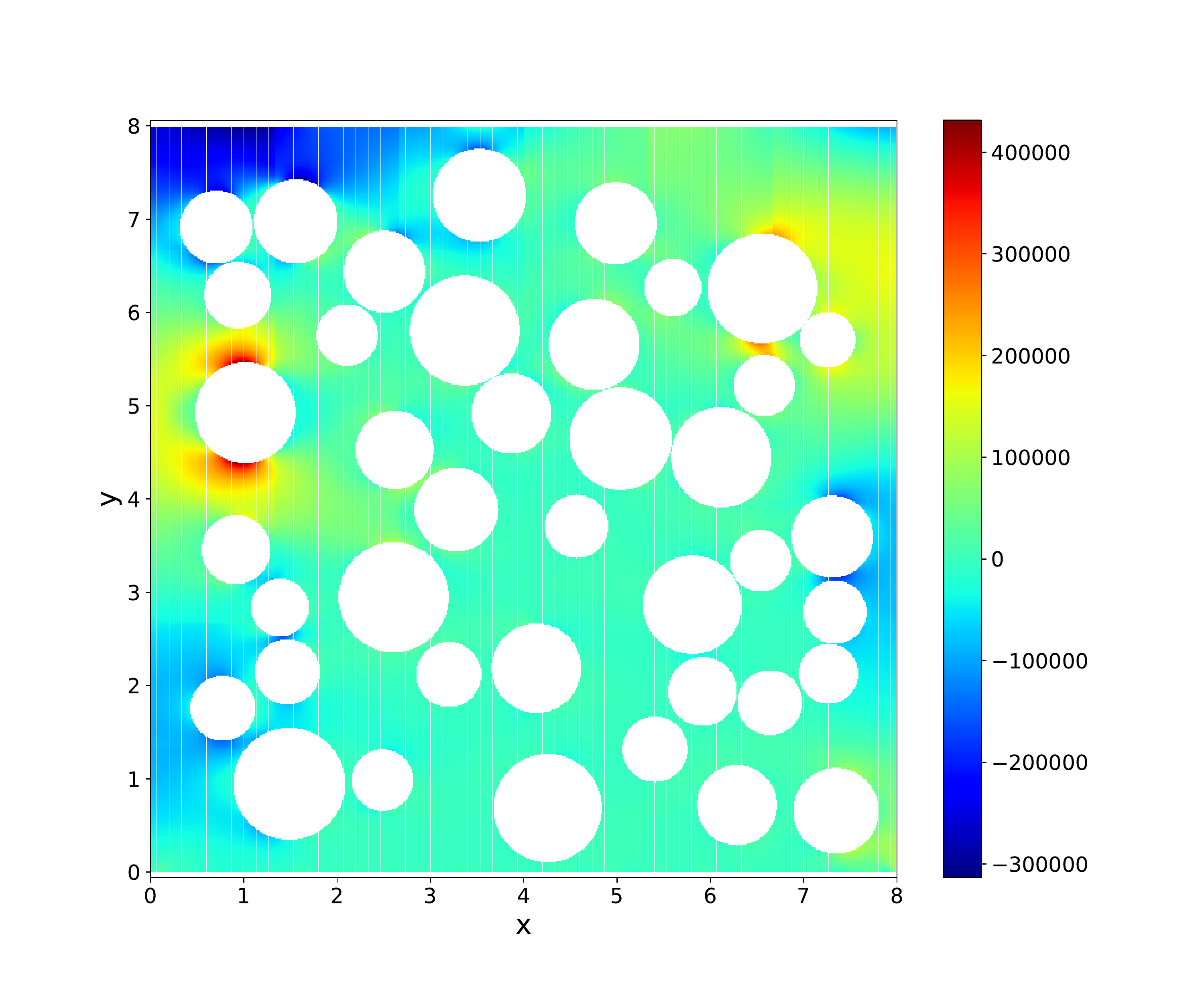}
		}
		\quad
		\subfigure[$\sigma_{x} \text{ over a cluster of $7$ holes}$]{
			\includegraphics[width=0.45\textwidth]{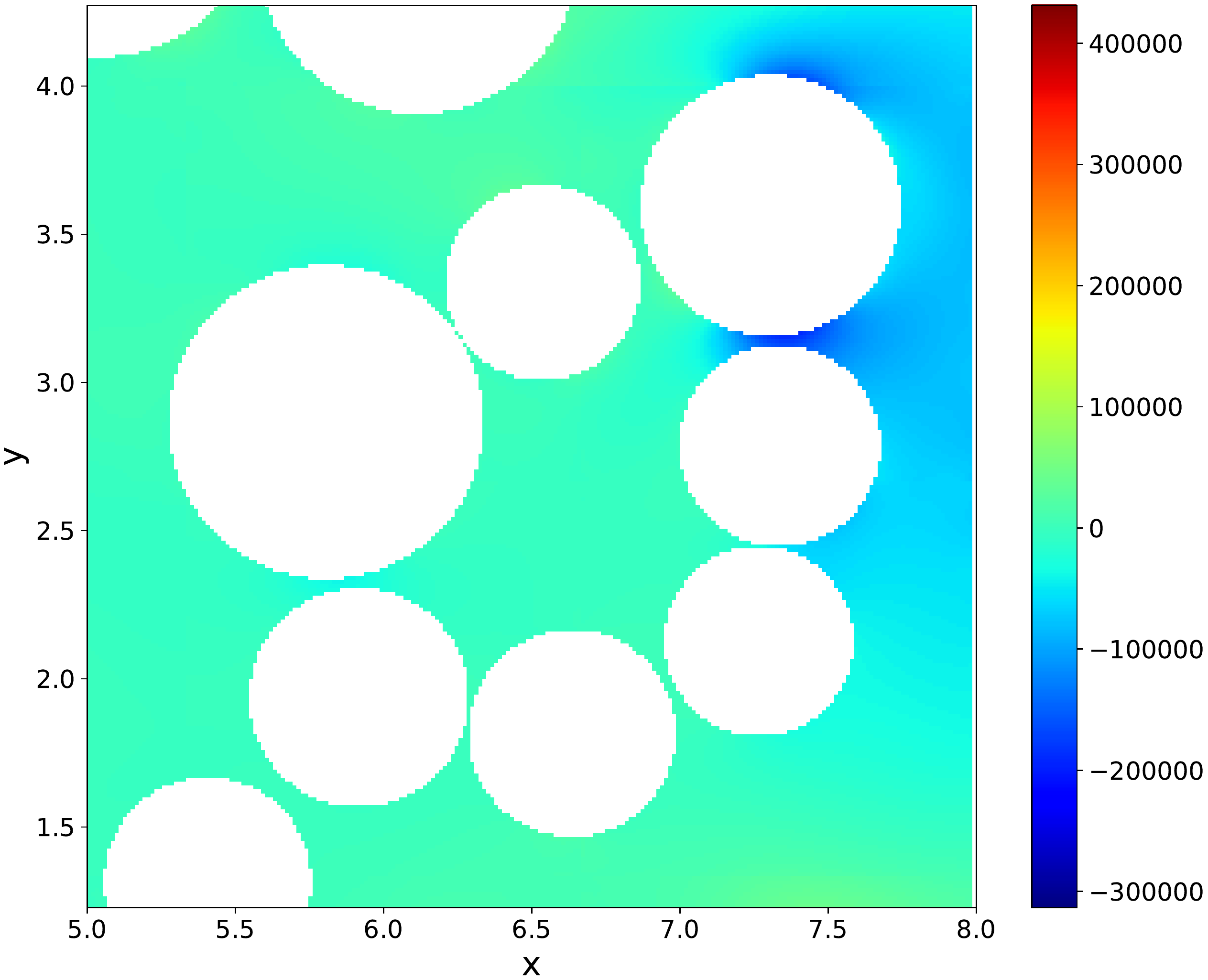}
		}
		\quad
		\subfigure[$\tau_{xy}$]{
			\includegraphics[width=0.45\textwidth]{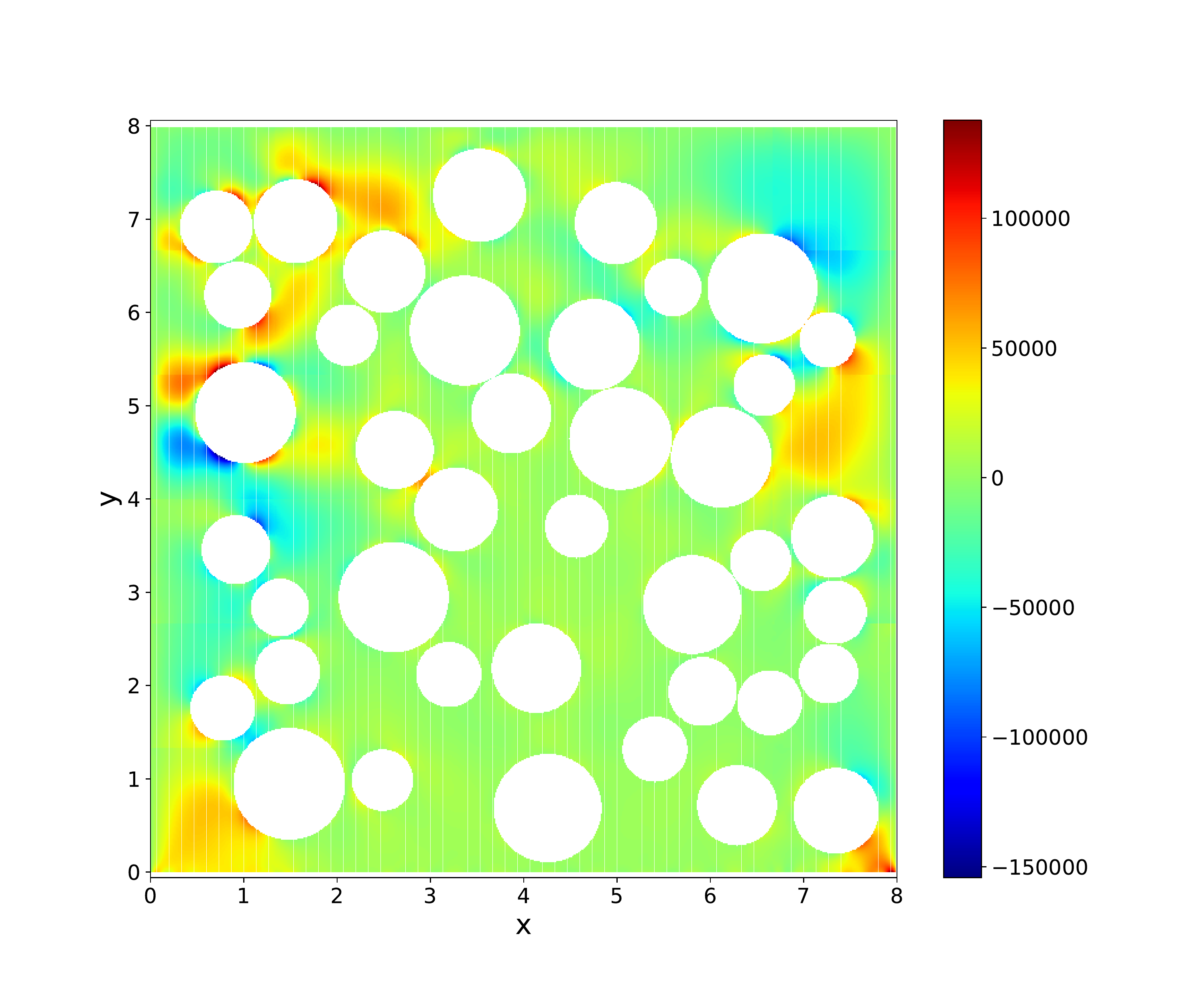}
		}
		\quad
		\subfigure[$\sigma_{y}$]{
			\includegraphics[width=0.45\textwidth]{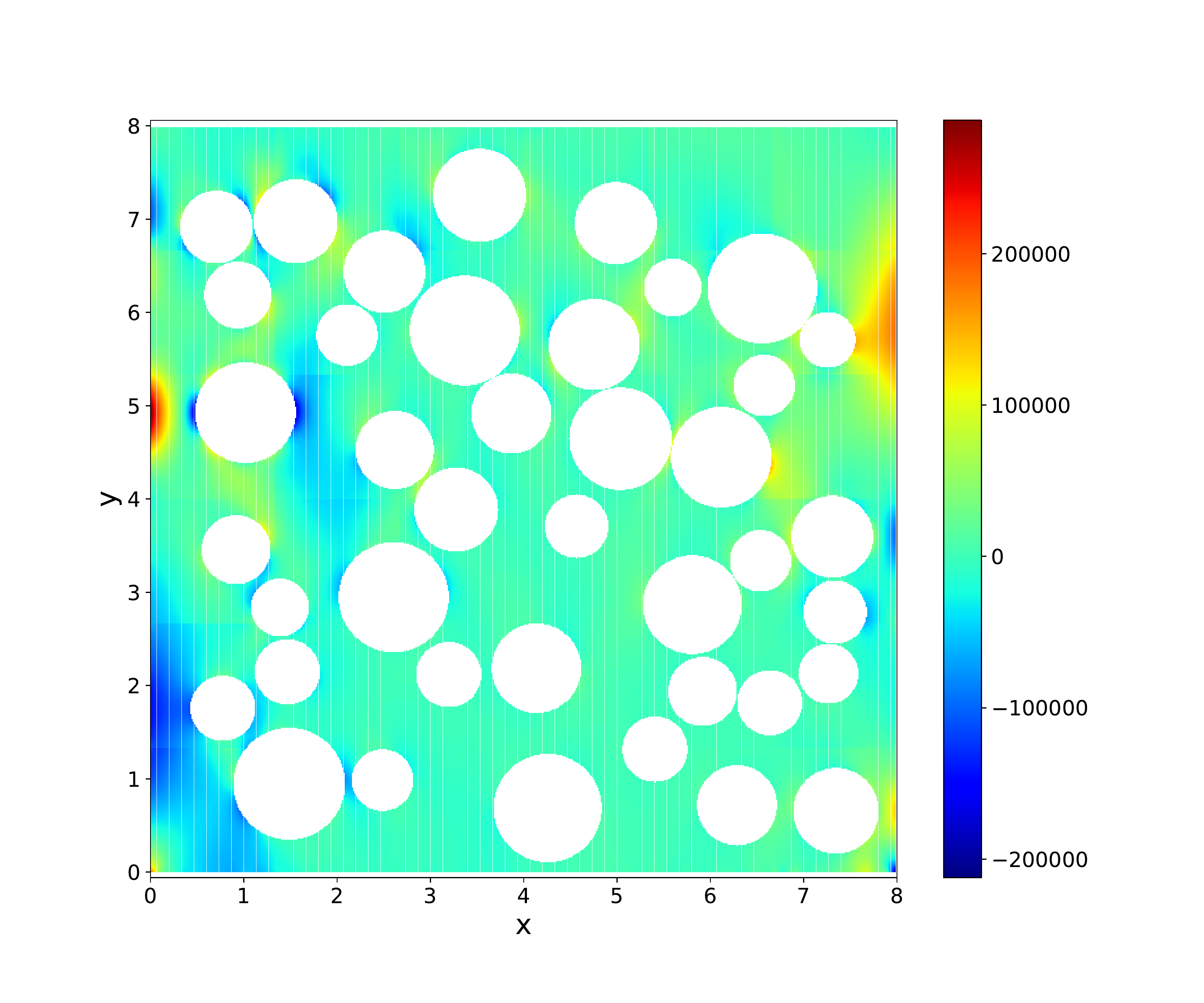}
		}
		\caption{Numerical solution by the random feature method for the two-dimensional elasticity problem over a complex geometry.}
		\label{fig3-3-2}
	\end{figure}
	Figure \ref{fig3-3-2} visualizes the displacement fields $u$, $v$, and the stress fields $\sigma_x$, $\tau_{xy}$, $\sigma_y$. 
	For this example,  it is quite difficult to generate a mesh for the FEM.
	If we simply remove the cluster in the lower right corner, we incur an $L^\infty$ error of about $50\%$ for $\sigma_{x}$; see Figure \ref{fig3-3-2}(c) and  \ref{fig3-3-2}(d).
	In contrast, it is quite straightforward to use  RFM to solve such a problem. As recorded in Table \ref{table3-3-2} in Appendix \ref{sec::appendix::complex::fem}, RFM shows a clear trend of numerical convergence.
	The error against the reference solution for displacements and stresses is reduced to about $5\%$. 
	
	\subsection{Homogenization}\label{sec3-4}
	
	In this subsection, we take a preliminary look at how RFM performs for problems with multi-scale solutions.
	We consider the elliptic equation over the unit disk 
	\begin{equation}
	\left\{\begin {aligned}
	-\operatorname{div}(a(\boldsymbol{x}) \nabla u(\boldsymbol{x}))  &= f(\boldsymbol{x}) \quad && \boldsymbol{x}\in\Omega, \\ 
	u(\boldsymbol{x})  &=  0 \quad && \boldsymbol{x}\in\partial \Omega,
	\end {aligned}\right.
	\end{equation}
	where $a(\boldsymbol{x})=e^{h(\boldsymbol{x})}$, $h(\boldsymbol{x})=\sum_{\mid \boldsymbol{k} \mid \leq R}(a_{\boldsymbol{k}} \sin (2 \pi \boldsymbol{k}  \cdot\boldsymbol{x})+b_{\boldsymbol{k}} \cos (2 \pi \boldsymbol{k} \cdot\boldsymbol{x}))$, $R=6$, and $\{a_{\boldsymbol{k}}\}$ and $\{b_{\boldsymbol{k}}\}$ are independent, identically distributed random variables with the distribution $\mathbb{U}[-0.3,0.3]$. 
	This is chosen so that there is no clear scale separation  in the coefficient \cite{owhadi2007metric}.
	
\begin{figure}[htbp]
		\centering
		\subfigure[$h$]{
			\includegraphics[width=0.3\textwidth]{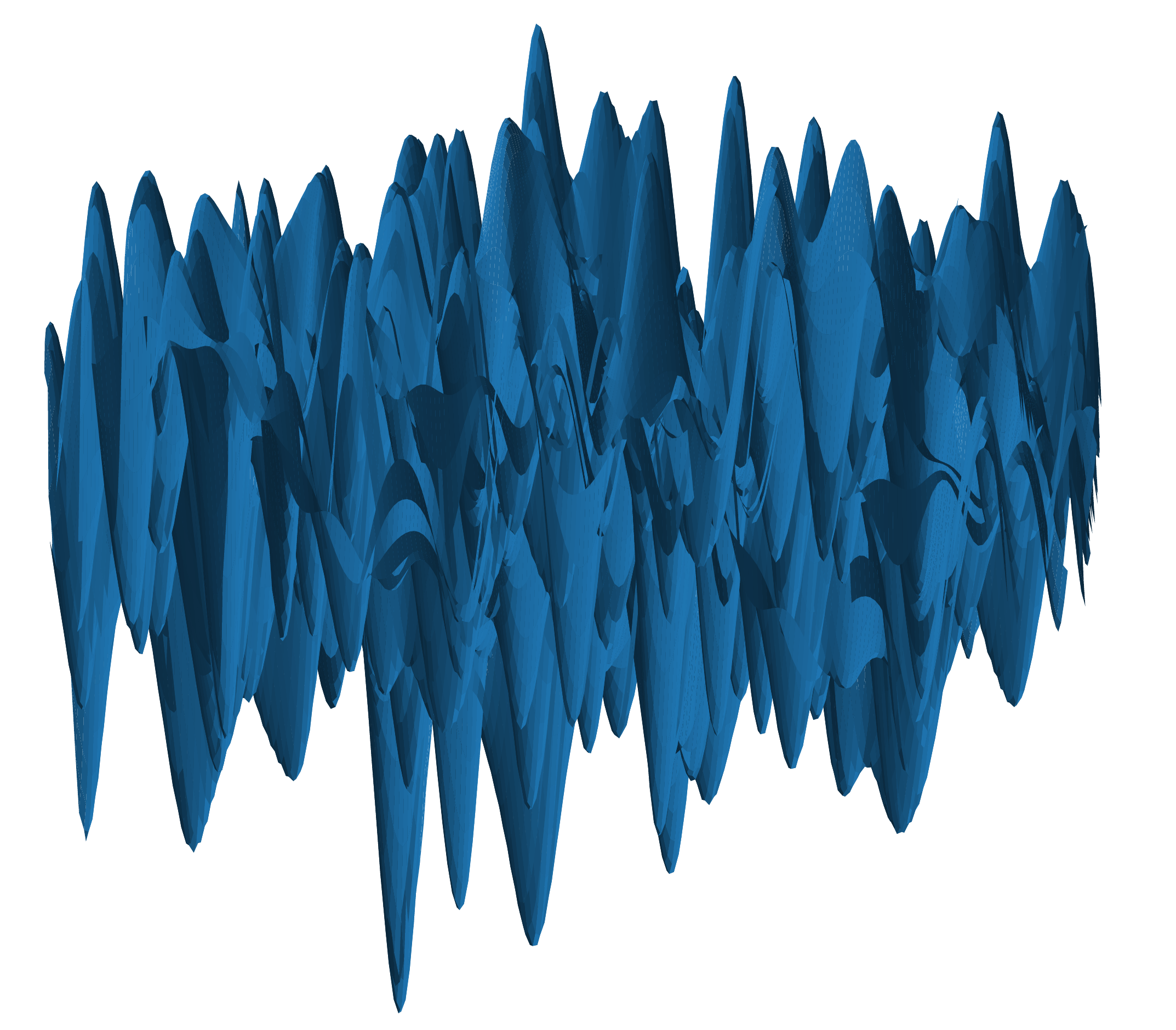}
		}
		\quad
		\subfigure[$a$]{
			\includegraphics[width=0.3\textwidth]{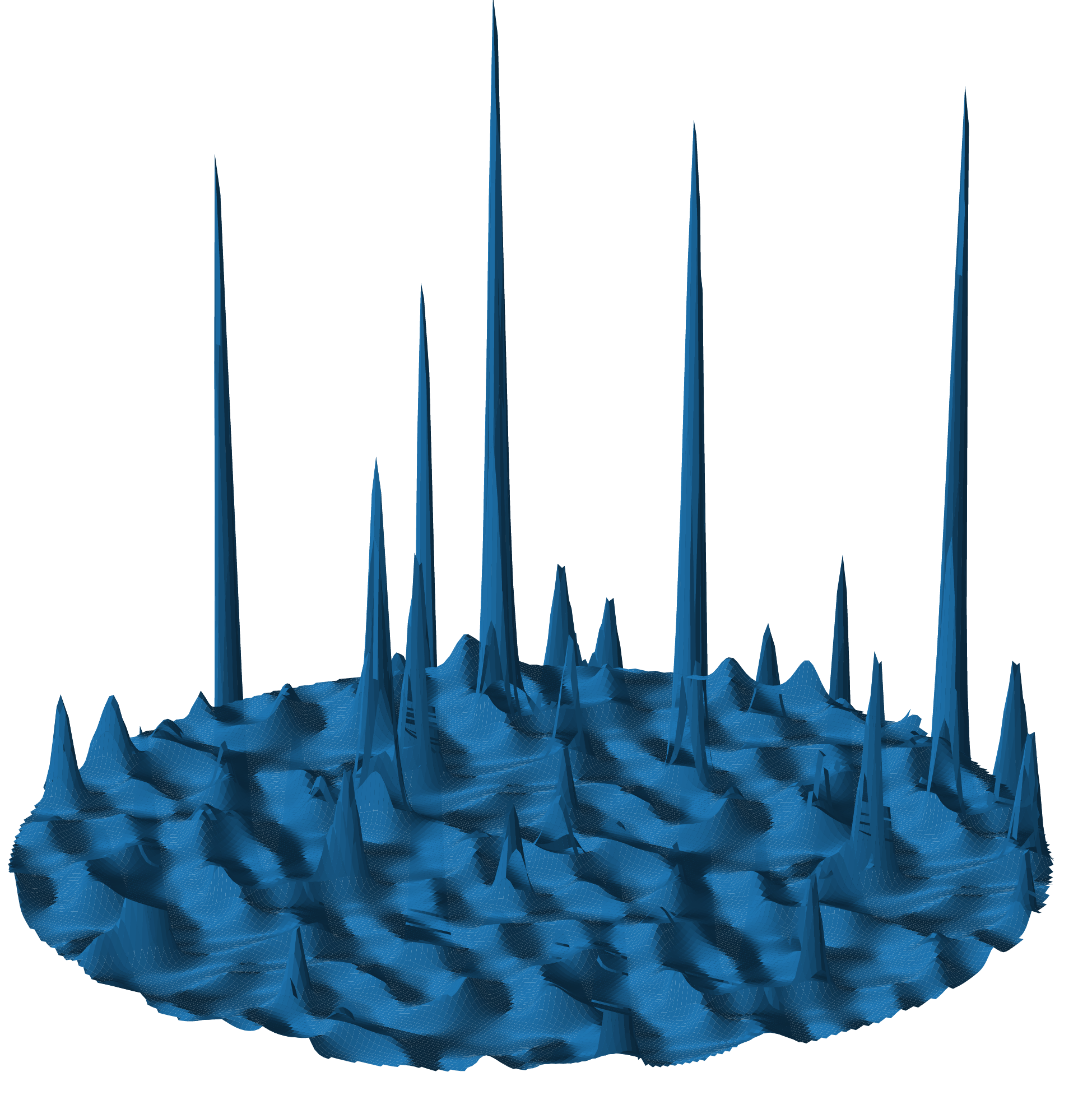}
		}
		\\
		\subfigure[$u$]{
			\includegraphics[width=0.28\textwidth]{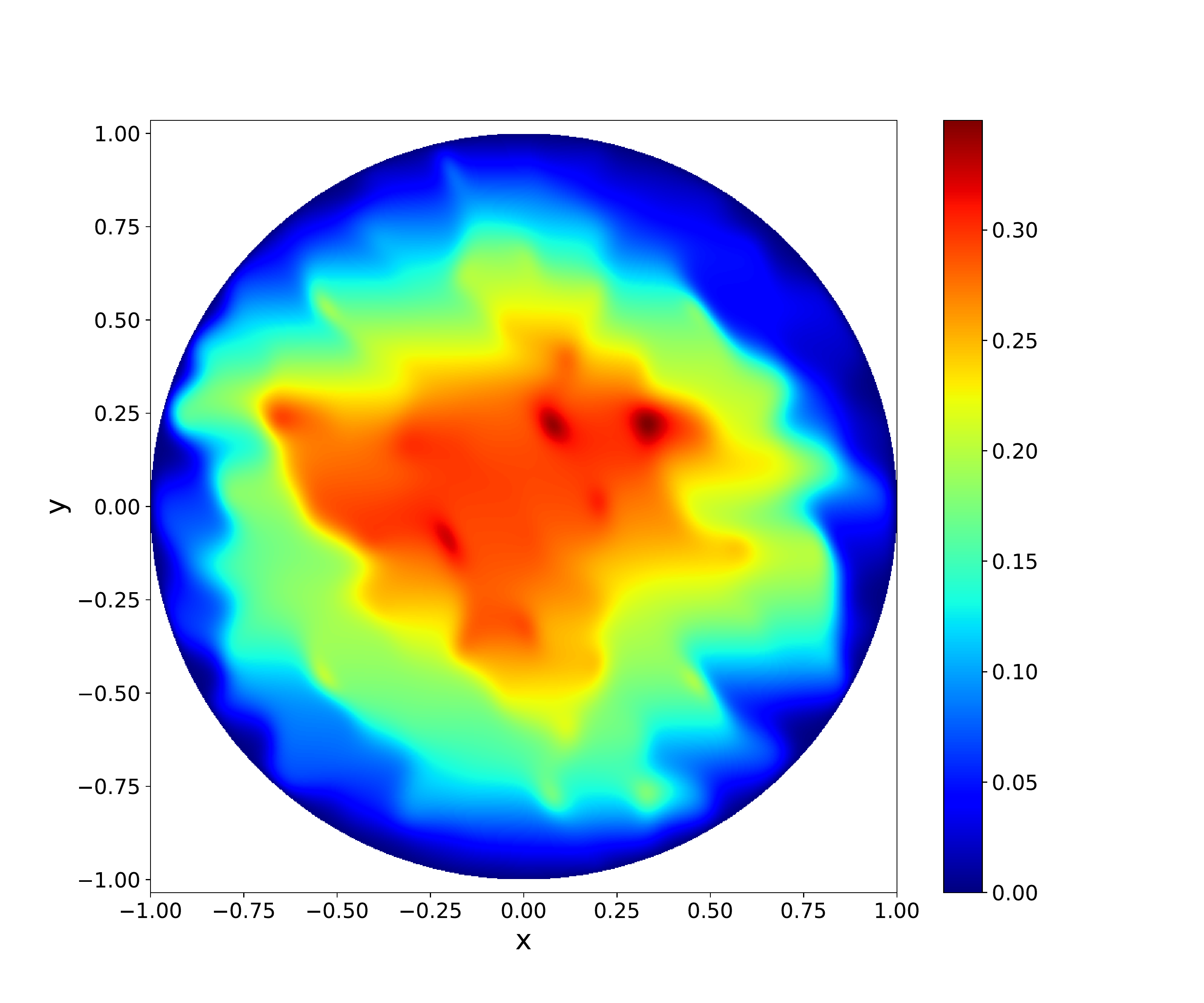}
		}
		\quad
		\subfigure[$u_{x}$]{
			\includegraphics[width=0.28\textwidth]{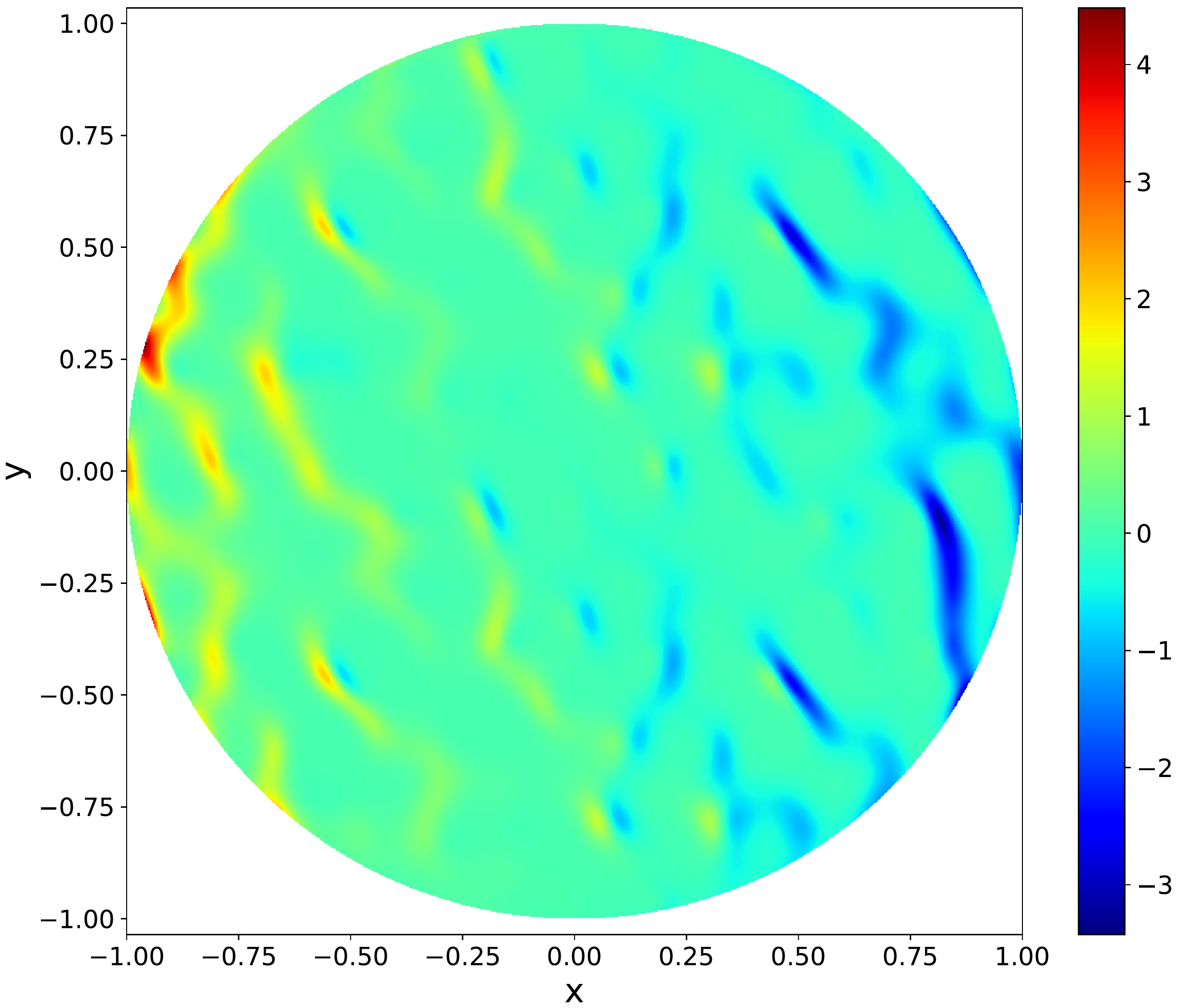}
		}
		\quad
		\subfigure[$u_{y}$]{
			\includegraphics[width=0.28\textwidth]{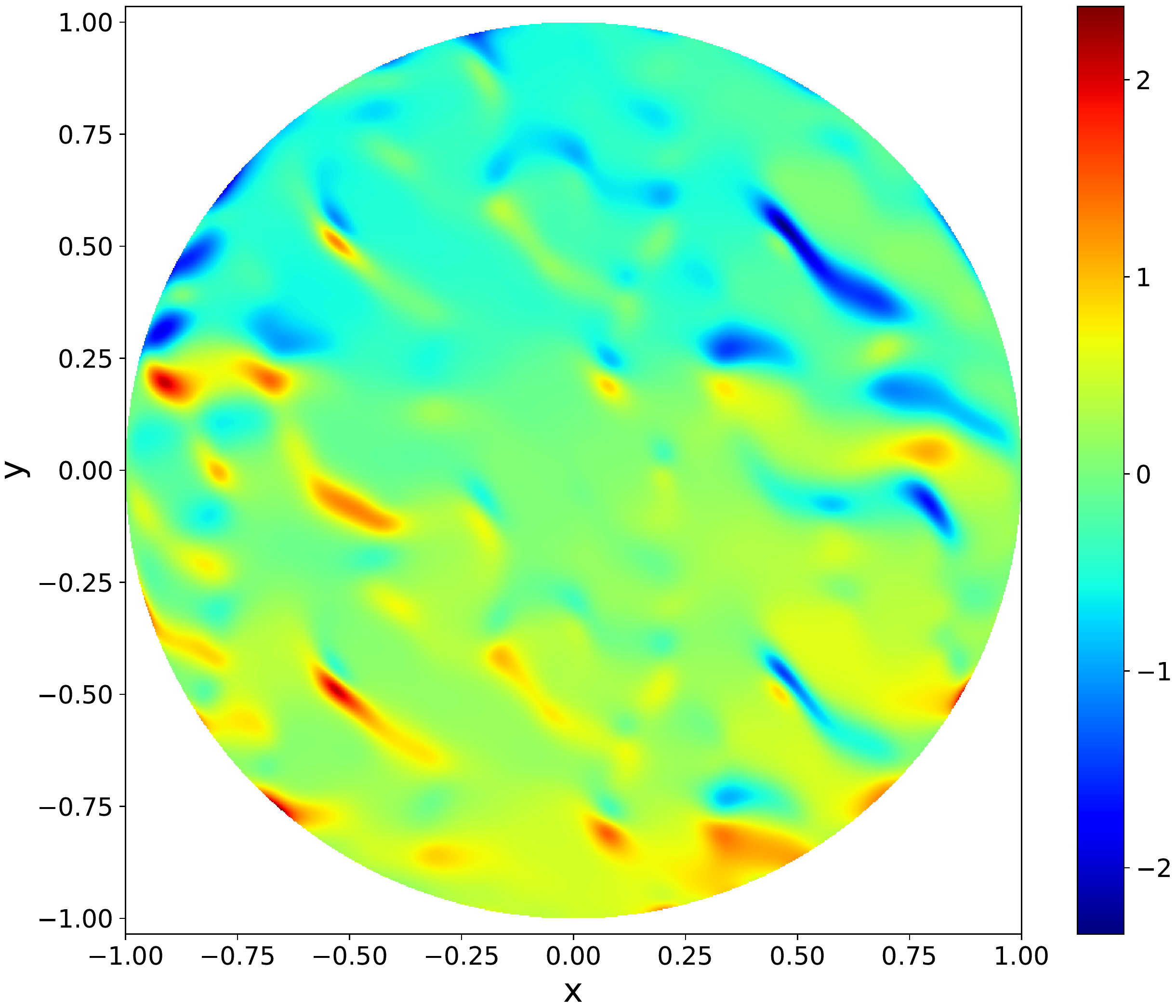}
		}
		\caption{Coefficient functions $h$ and $a$, the numerical solution and its first-order derivatives obtained by the random feature method for the homogenization problem.}
		\label{fig3-4}
	\end{figure}

	More details of the experimental setup can be found in Appendix \ref{sec::appendix::complex::homogenization}. Figure \ref{fig3-4} visualizes the coefficient functions $h$ and $a$, the numerical solution and its first-order derivatives obtained by RFM. Table \ref{table3-4} in Appendix \ref{sec::appendix::complex::homogenization} records the convergence 
	behavior of RFM when the solution with $N=86219$ is taken as the reference. 

	\subsection{Stokes flow}\label{sec3-5}
	Consider  Stokes flow  defined by 
	\begin{equation}
	\left\{\begin {aligned}
	-\Delta \boldsymbol{u}(\boldsymbol{x}) + \nabla p(\boldsymbol{x})&=\boldsymbol{f}(\boldsymbol{x}) \quad && \boldsymbol{x}\in\Omega, \\ 
	\nabla\cdot \boldsymbol{u}(\boldsymbol{x}) &=0  \quad && \boldsymbol{x}\in \Omega,\\
	\boldsymbol{u}(\boldsymbol{x})  &= \boldsymbol{U}(\boldsymbol{x}) \quad && \boldsymbol{x} \in \partial\Omega.
	\end {aligned}\right.
	\label{stokesflow}
	\end{equation}
	In this problem,  $p$ is only determined up to a constant. To avoid difficulties, we fix the value of $p$ at the left-bottom corner.
	
	One problem with spectral methods is that 
	spurious pressure mode arises due to the rank deficiency of the discrete systems  \cite{schumack1991spectral}. 
	One interesting feature of RFM is that it always looks for an optimal solution with minimal norm. 
	This allows us to automatically bypass the issue of rank deficiency, as we see in the following examples. 
	
	First, we consider \eqref{stokesflow} with an explicit solution and inhomogeneous boundary condition, where $\Omega$ is 
	the square $(0,1) \times(0,1)$ with three holes centered at $(0.5,0.2)$, $(0.2,0.8)$, $(0.8,0.8)$ of radius $0.1$. The exact displacement fields and the experimental setup are detailed  in Appendix \ref{sec::appendix::complex::stokes}. Table \ref{table3-5-1} in Appendix \ref{sec::appendix::complex::stokes} records the convergence behavior of RFM and spectral accuracy is observed for $u$, $v$ as well as $p$.
	
\begin{figure}[htbp]
		\centering
		\subfigure[]{
			\includegraphics[width=0.40\textwidth]{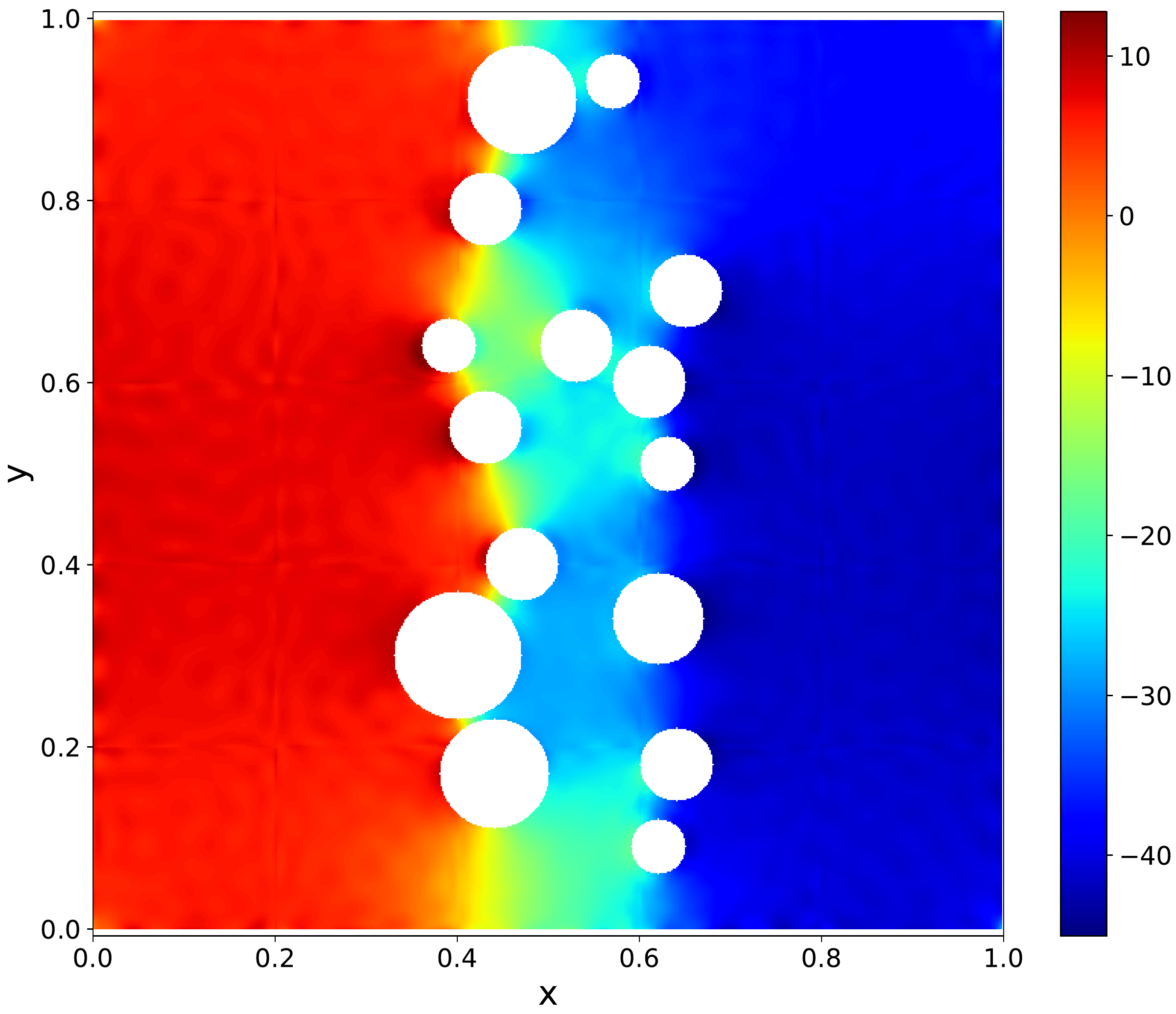}
		}
		\quad
		\subfigure[]{
			\includegraphics[width=0.40\textwidth]{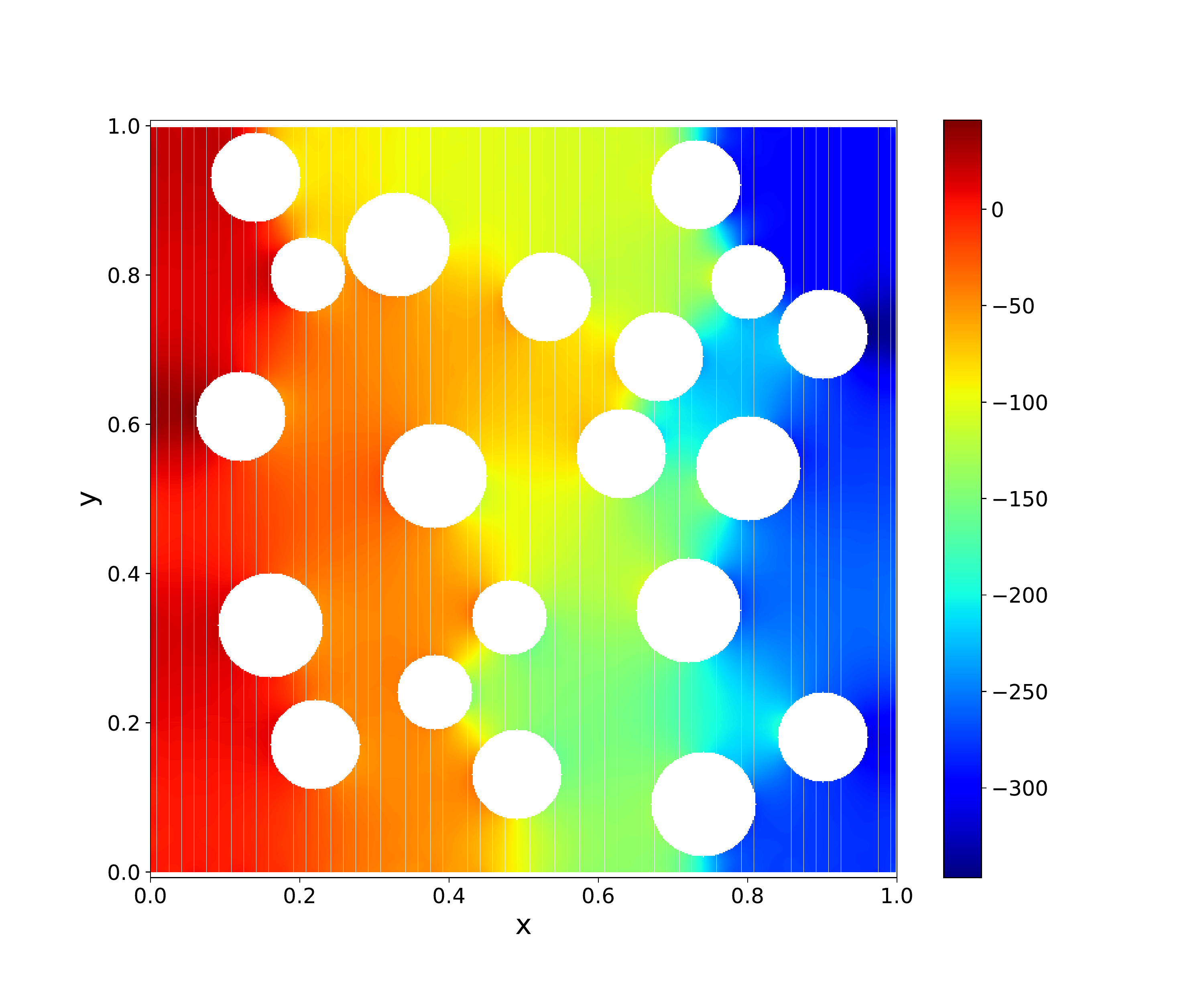}
		}
		\quad
		\subfigure[]{
			\includegraphics[width=0.40\textwidth]{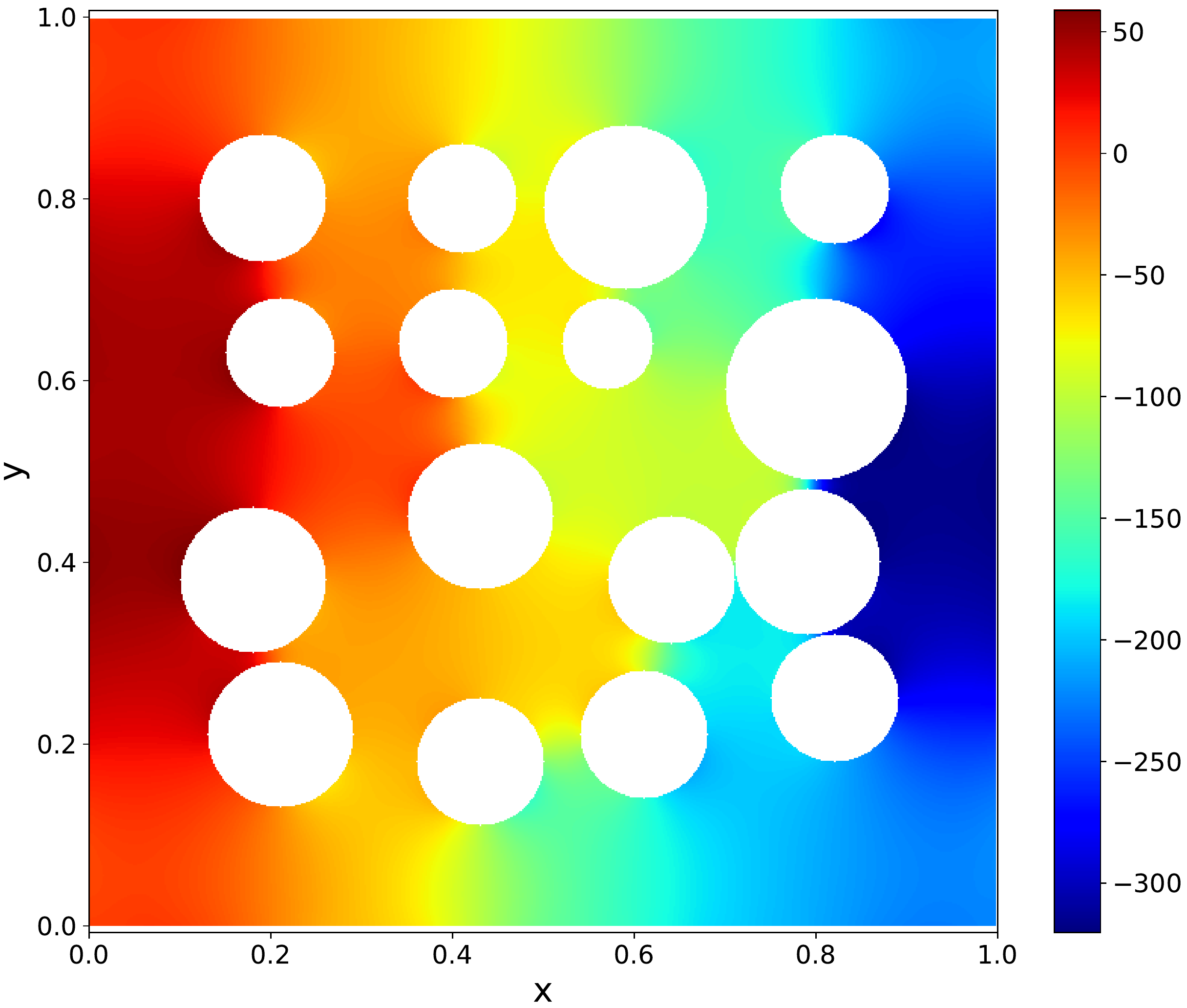}
		}
		\quad
		\subfigure[]{
			\includegraphics[width=0.40\textwidth]{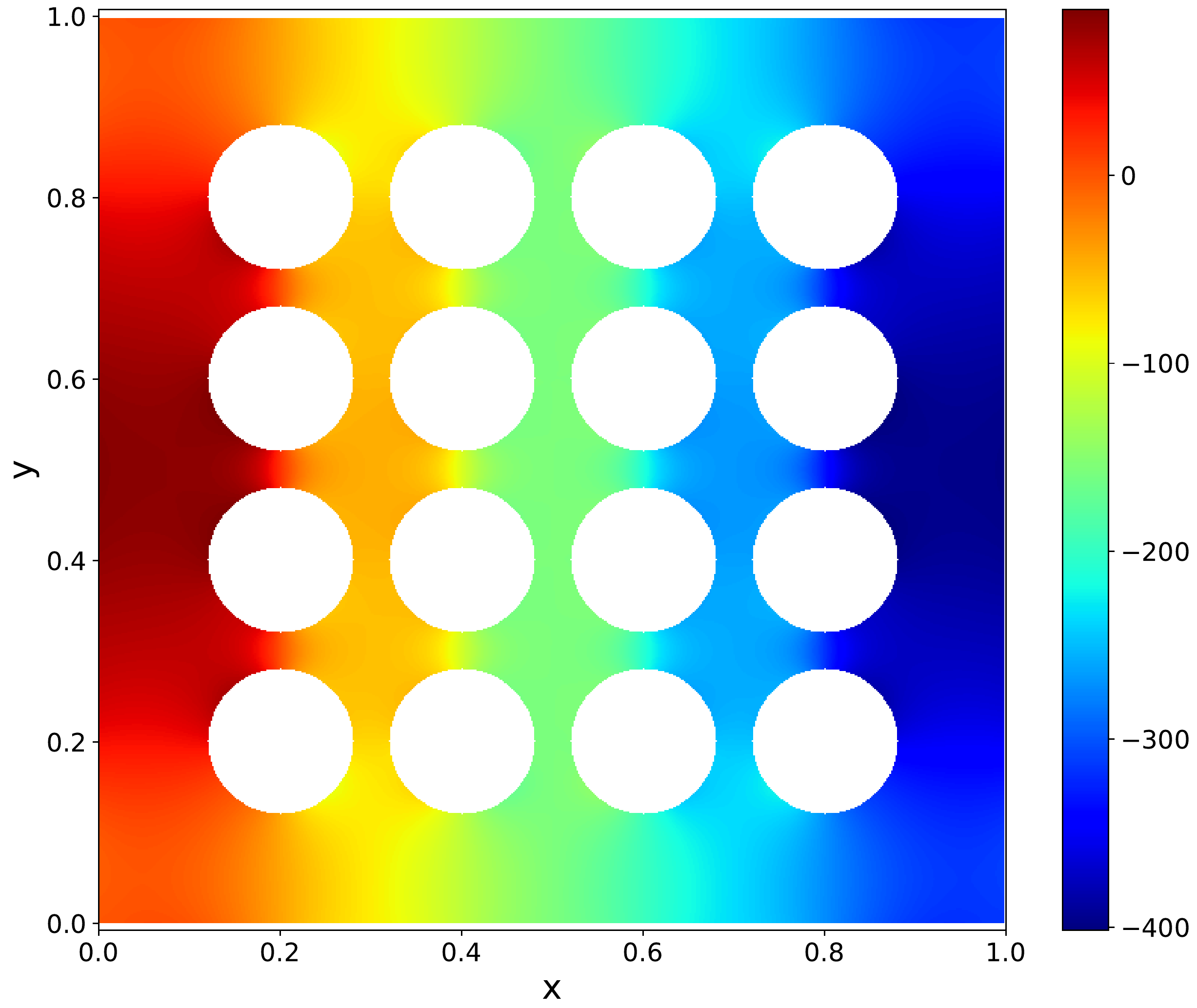}
		}
		\caption{Pressure diagram generated by the random feature method for four sets of complex obstacles.}\label{fig3-5-1}
	\end{figure}

	Next, we consider two-dimensional channel flows for four sets of complex obstacles with the inhomogeneous boundary condition
	\begin{equation}
	(u,v)|_{\partial \Omega} = 
	\begin{cases}
	(y  (1-y),0) & \text { if } x=0 \\ 
	(y  (1-y),0) & \text { if } x=1 \\ 
	(0,0) & \text {otherwise}
	\end{cases}.
	\end{equation}
	The pressure diagram is plotted in Figure \ref{fig3-5-1}.
	\begin{figure}[htbp]
		\centering
		\subfigure[$u$]{
			\includegraphics[width=0.45\textwidth]{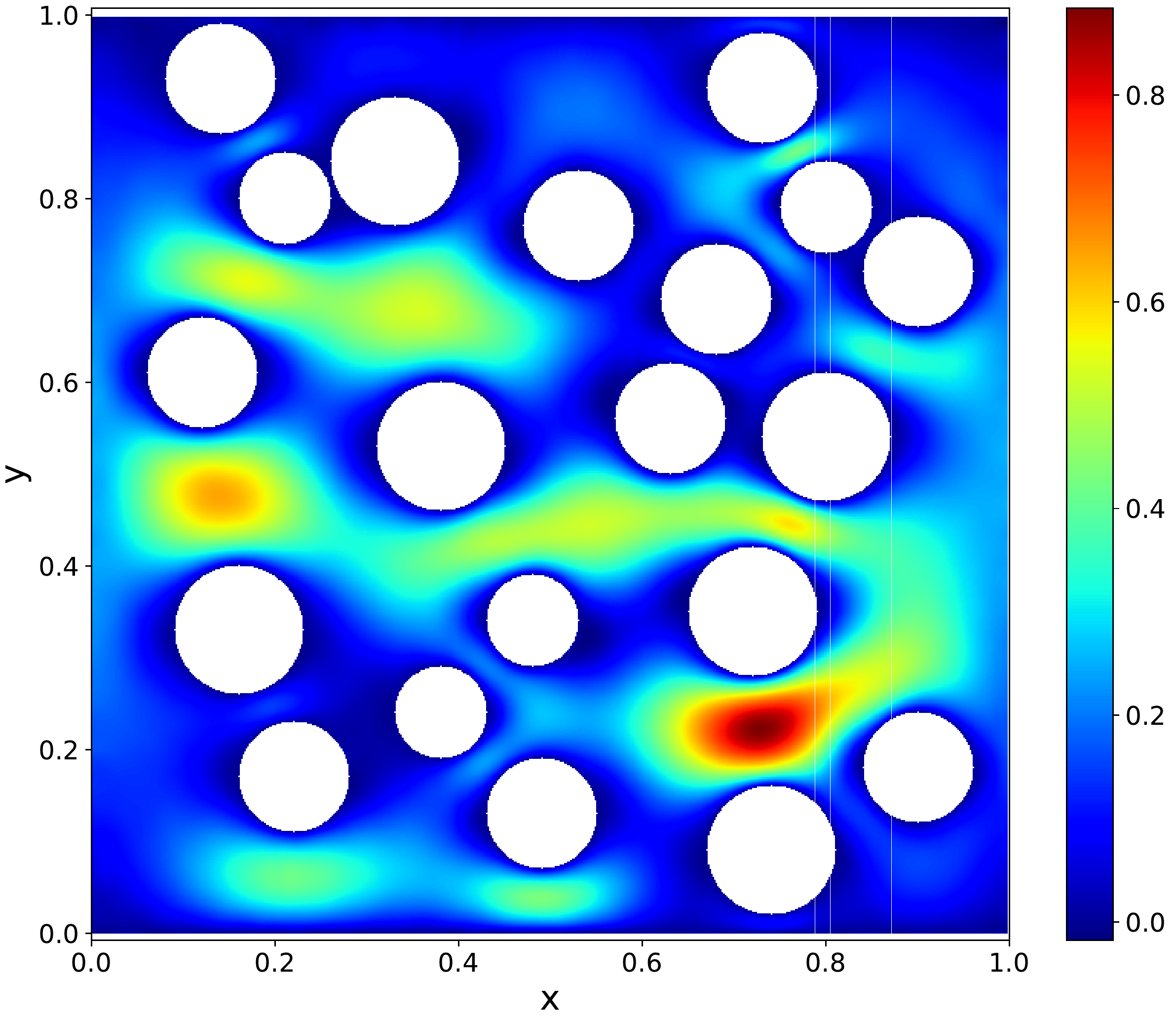}
		}
		\quad
		\subfigure[$v$]{
			\includegraphics[width=0.45\textwidth]{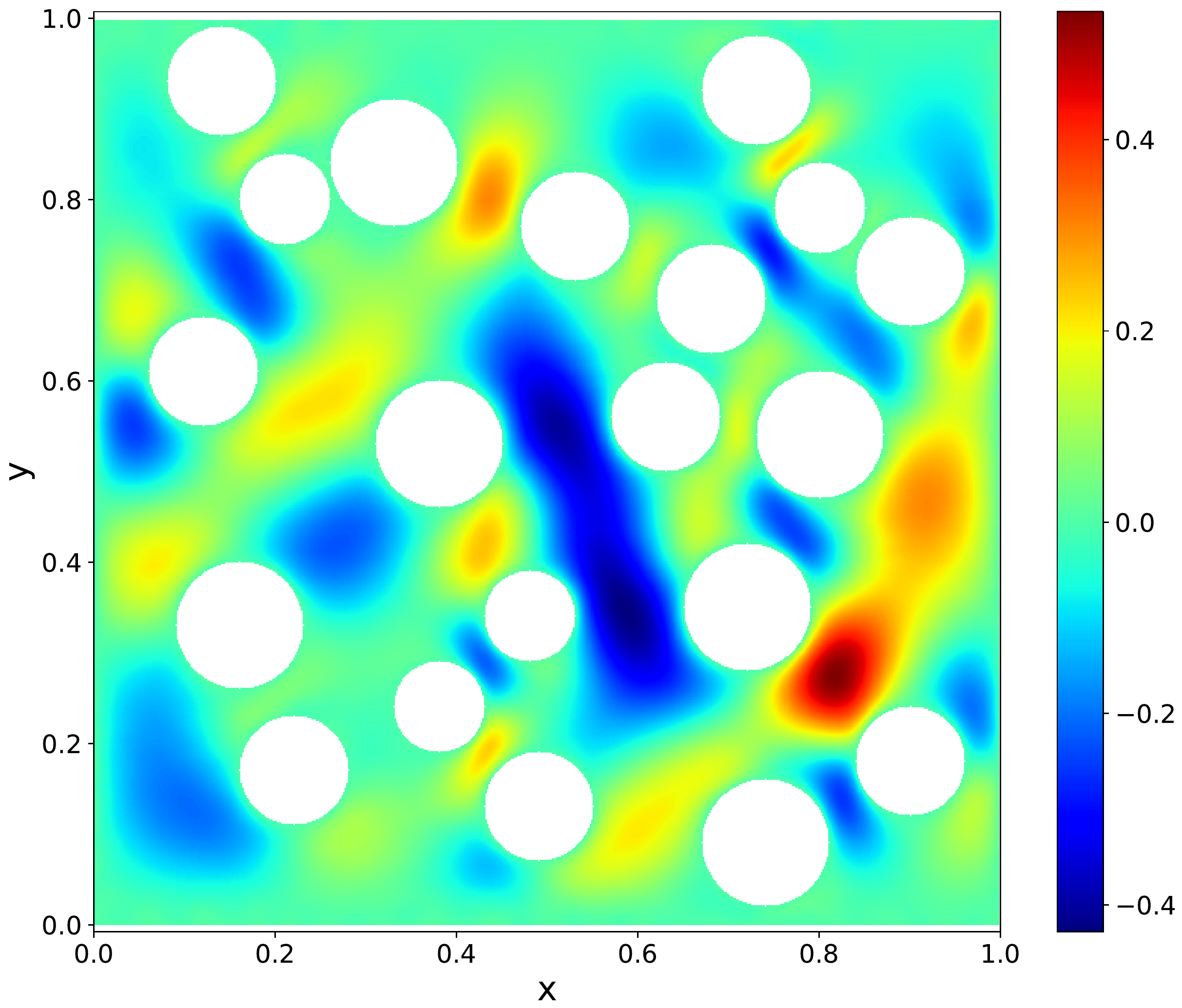}
		}
		\caption{Velocity field $(u, v)$ generated by the random feature method.}\label{fig3-5-2}
	\end{figure}
	Numerical convergence is observed for all four examples and the difference with respect to the reference solution is about $1\%$ for $u$, $v$, $p$. The difference is about $0.01\%$ for uniformly distributed holes.
	
	\section{Discussions}\label{sec4}
	
	One main disadvantage of the traditional algorithms is that they are not flexible enough.
	This lack of flexibility is reflected in several aspects. For example, most traditional algorithms
	require that the number of free parameters be the same as the number of conditions.
	Spectral methods often require that the basis functions satisfy particular boundary conditions,
	and that the basis functions are constructed in a tensor-product form.
	
	Neural network-based methods such as  Deep BSDE \cite{EHJ2017}, Deep Ritz \cite {EY2018} and PINN \cite{PINN} do not have these problems.  However, they typically lack robustness.
	With these methods, it is often quite easy to get a roughly accurate solution, but very hard to systematically improve the accuracy.
	As a result,  one is only willing to use them  when traditional algorithms fail to work.
	
	The random feature method introduced here seems to have both the flexibility and robustness needed.
	There are two major differences of RFM compared with neural network-based methods such as Deep BSDE, Deep Ritz or PINN. The first is that we use only the random feature model  instead of the neural network model to represent the solution. 
	We find that for these low dimensional problems, there is not much gain in terms of the representative power by 
	switching to neural networks,
	yet the non-convexity introduced in the loss function makes the task of training much harder.
	The second is that we use a multi-scale representation instead of a single global representation.
	Compared with traditional algorithms, an important difference is that RFM typically works in the situation where the number of unknown parameters
	(say $M$) is different from the number of conditions (say $N$). This forces us to use a least square framework, which increases the complexity
	of the training process since the condition number is now much bigger.
	In return, it allows us to obtain reasonable solutions with much lower human and computer cost. For example, in the complex geometry problems
	treated earlier, we need a large value of $N$ to resolve the geometry.
	If we used the same value of $M$,  the cost would be too big in our current implementation due to the usage of \textit{scipy} in Python.
	
	
	Another important difference compared with traditional algorithms is that the boundary condition is treated in the same fashion as the PDE. 
	In particular, we do not force the basis functions to satisfy particular boundary conditions.
	This increased flexibility is vital for the success on problems with complex geometry.
	
	A third important difference is the adoption of neural network-like basis functions instead of traditional tensor-product based basis functions.
	This means that the number of terms needed does not necessarily go up like $n^d$ with $n$ being the number of unknown parameters in each dimension and $d$ being the dimensionality.
	Instead it depends  entirely on the complexity of the solution. This is also part of the reason why random choices of the feature vectors
	is generally preferred.
	
	RFM differs from the local extreme learning machines in the following aspects. 
	The first is that a partition of unity is used to construct the local random feature functions instead of domain decomposition. 
	This allows us to bypass the smoothness conditions and thus simplifies the loss function and the subsequent training. 
	The second is the rescaling procedure. Though the difference for model problems is small, we find it crucial for
	practical problems of interest such as the elasticity problem or the Stokes flow problem.
	
	Let us now turn to a discussion of the crucial components in RFM. The first is the choice of the basis functions.
	Here a crucial issue is the probability distribution for the feature vector.
	We can use pre-computing to give us a rough idea about this distribution.  
	In practice, we find that as long as the support of the distribution approximately covers the  frequency range of the
	true solution, RFM produces stable results with sine/cosine activation functions. 
	In addition, we observe that in most cases, random sampling performs better than deterministic choices of the feature vectors.
	For example, results for the Helmholtz equation in Section \ref{sec3-1-3} show that even when $d = 1$, in most cases,  random
	sampling of $k$ and $b$ performs better than choosing the values of $k$ and $b$ from a uniform grid, and is at least
	as good in the remaining cases.
	This might be counterintuitive to what we have learned in classical numerical analysis that quadrature schemes are superior to the Monte-Carlo method for low-dimensional integrals (say $d\leq 3$). 
	We are in the process of trying to quantify this finding.
	
	The second is the choice of collocation points. 
	Ideally we would like the collocation points to be equally distributed, both in the interior and at the boundary.
	This becomes non-trivial for three dimensional situations when the boundary is a surface.
	We are in the process of developing techniques that can help us to accomplish this.
	
	The third technical aspect is the training.  By turning to a least square formulation we may have increased the size of the
	condition number. There are a number of preconditioning and reformulation techniques that might be useful to alleviate this problem.
	In any case, it would be useful to carry out a precise numerical analysis of carefully chosen model problems to gain some insight about
	the convergence behavior of the training process.
	
	
	\section*{Acknowledgments}
	We thank Prof. Suchuan Dong for providing the one-dimensional code of local extreme learning machines and for helpful discussion. The work is supported by Anhui Center for Applied Mathematics, and the Major Project of Science \& Technology of Anhui Province (No. 202203a05020050). J. Chen also acknowledges supported by NSFC 11971021.

	\begin{appendix}
		\section{Numerical results for different choices of random feature functions}
		\label{sec::appendix::basis}
		\subsection{Partition of unity and local random feature models}
		\label{sec3-1-1}
		
		The introduction of PoU generates local random feature functions and provides a more general strategy than domain decomposition and mesh generation.
		
		Consider the following explicit solution to \eqref{helmholtz}
		\begin{equation}
		u(x)=\sin (3 \pi x+\frac{3 \pi}{20}) \cos (2 \pi x+\frac{\pi}{10})+2.
		\label{helmholtz-sol1}
		\end{equation}
		In a series of tests, we set the hyper-parameters as follows:
		\begin{itemize}
			\item $M=200$, $400$, $800$, $1600$;
			
			\item $J_{n}=50$;
			
			\item $M_p=\frac{M}{50}$;
			
			\item $\boldsymbol{x}_{n} = x_n = 8\frac{2n-1}{2M_p}, n=1,\cdots,M_p$;
			
			\item $\boldsymbol{r}_{n} = r_{n}=\frac{8}{2M_p}, n=1,\cdots,M_p$;
			
			\item $Q=200$, $400$, $800$, $1600$.
		\end{itemize}
		For $\psi^{a}$, we need additional $2(M_p-1)$ smoothness conditions for the solution and its first derivative at $x=\frac{1}{M_p},\cdots,\frac{M_p-1}{M_p}$. No additional conditions are needed for $\psi^{b}$. The total number of conditions is $N=208$, $416$, $832$, $1664$ for $\psi^{a}$,  and $N=202$, $402$, $802$, $1602$ for $\psi^{b}$, respectively.
		
		Table \ref{table3-1-1-1} compares RFM and PINN \cite{PINN} in terms of accuracy. 
		The network used in PINN has the same structure as that in RFM, one hidden layer with $M$ neurons and the $\tanh$ activation function is used. Collocation points are also chosen to be the same as those in RFM with $\psi^{b}$. Since the inner parameters are trainable in PINN, we use the Adam optimizer with  learning rate $0.001$ to train the network. The training process ends after $100000$ epochs when $M=200$, $400$ and $200000$ epochs when $M=800$, $1600$.
		
		\begin{table}[htbp]
			\caption{Comparison of the RFM and PINN  for the one-dimensional Helmholtz equation.}\label{table3-1-1-1}
			\centering
			\begin{small}
				\begin{tabular}{|c|cc|cc|cc|}
					\hline
					\multirow{2}{*}{M} & \multicolumn{2}{|c|}{$\psi^{a}$} &  \multicolumn{2}{|c|}{$\psi^{b}$} & \multicolumn{2}{|c|}{PINN} \\
					\cline{2-7}
					& N & $L^{\infty}$ error & N & $L^{\infty}$ error & N & $L^{\infty}$ error \\
					\hline
					200  & 208  & 8.76E-2  & 202  & 2.51E-2  & 202  & 2.59E-2  \\
					400  & 416  & 5.89E-7  & 402  & 5.18E-7  & 402  & 6.77E-3  \\
					800  & 832  & 4.44E-10 & 802  & 6.61E-10 & 802  & 1.35E-2  \\
					1600 & 1664 & 8.84E-12 & 1602 & 1.18E-11 & 1602 & 8.94E-3  \\
					\hline
				\end{tabular}
			\end{small}
		\end{table}
		
		From Table \ref{table3-1-1-1}, we observe that error in  PINN is around $1E-3$ without notable further improvement, while RFM for different PoU functions has exponential convergence. 
		This suggests that fixing the inner parameters greatly simplifies the optimization problem and allows us to obtain accurate and robust solutions.
		
		Next, we report the results for Poisson equation \eqref{poisson} with the following explicit solution
		\begin{equation}
		u(x, y)=-[\frac{3}{2} \cos (\pi x+\frac{2\pi}{5})+2 \cos (2 \pi x-\frac{\pi}{5})][\frac{3}{2} \cos (\pi y+\frac{2\pi}{5})+2 \cos (2 \pi y-\frac{\pi}{5})].
		\label{poisson-sol1}
		\end{equation}
		
		In this problem, we set the  hyper-parameters as follows:
		\begin{itemize}
			\item $M=200$, $400$, $800$, $1600$;
			
			\item $J_{n}=400$;
			
			\item $M_p=\frac{M}{400}$;
			
			\item $\boldsymbol{x}_{n} = (x_i, y_j) = (\frac{2i-1}{2\sqrt{M_p}},\frac{2j-1}{2\sqrt{M_p}}), i,j=1,\cdots,\sqrt{M_p}$;
			
			\item $\boldsymbol{r}_{n} = r_{i,j}=(\frac{1}{2\sqrt{M_p}},\frac{1}{2\sqrt{M_p}}), i,j=1,\cdots,\sqrt{M_p}$;
			
			\item $Q=400M_p$, $625M_p$, $900M_p$, $1225M_p$, $1600M_p$.
		\end{itemize}
		For $\psi^{a}$, we impose the smoothness conditions at $(\frac{i}{\sqrt{M_p}}, \frac{j}{\sqrt{M_p}}), i,j=1,2,\cdots,\sqrt{M_p}$. 
		Table \ref{table3-1-1-2} shows the error of RFM for different PoU functions. Both show exponential convergence.
		
		\begin{table}[htbp]
			\caption{Results  for the two-dimensional Poisson equation with explicit solution \eqref{poisson-sol1}.}\label{table3-1-1-2}
			\centering
			\begin{small}
				\begin{tabular}{|c|cc|cc|}
					\hline
					\multirow{2}{*}{M} & \multicolumn{2}{|c|}{$\psi^{a}$} & \multicolumn{2}{|c|}{$\psi^{b}$} \\
					\cline{2-5}
					& N & $L^{\infty}$ error & N & $L^{\infty}$ error \\
					
					\hline
					\multirow{5}{*}{1600} 
					& 1920 & 1.74E-8  & 1760 & 1.90E-7  \\
					& 2900 & 1.55E-9  & 2700 & 1.22E-10 \\
					& 4080 & 2.31E-10 & 3840 & 3.89E-10 \\
					& 5460 & 6.29E-11 & 5180 & 2.67E-10 \\
					& 7040 & 5.04E-11 & 6720 & 4.68E-10 \\
					
					\hline
					\multirow{5}{*}{6400} 
					& 7680  & 7.77E-10 & 6720  & 1.61E-8   \\
					& 11600 & 5.74E-11 & 10400 & 1.91E-11  \\
					& 16320 & 7.04E-12 & 14880 & 5.64E-11  \\
					& 21840 & 9.93E-12 & 20160 & 5.21E-11  \\
					& 28160 & 1.66E-11 & 26240 & 4.97E-11  \\
					\hline
				\end{tabular}
			\end{small}
		\end{table}
		
\begin{figure}[htbp]
			\centering
			\subfigure[$\psi^{a}$]{
				\includegraphics[width=0.45\textwidth]{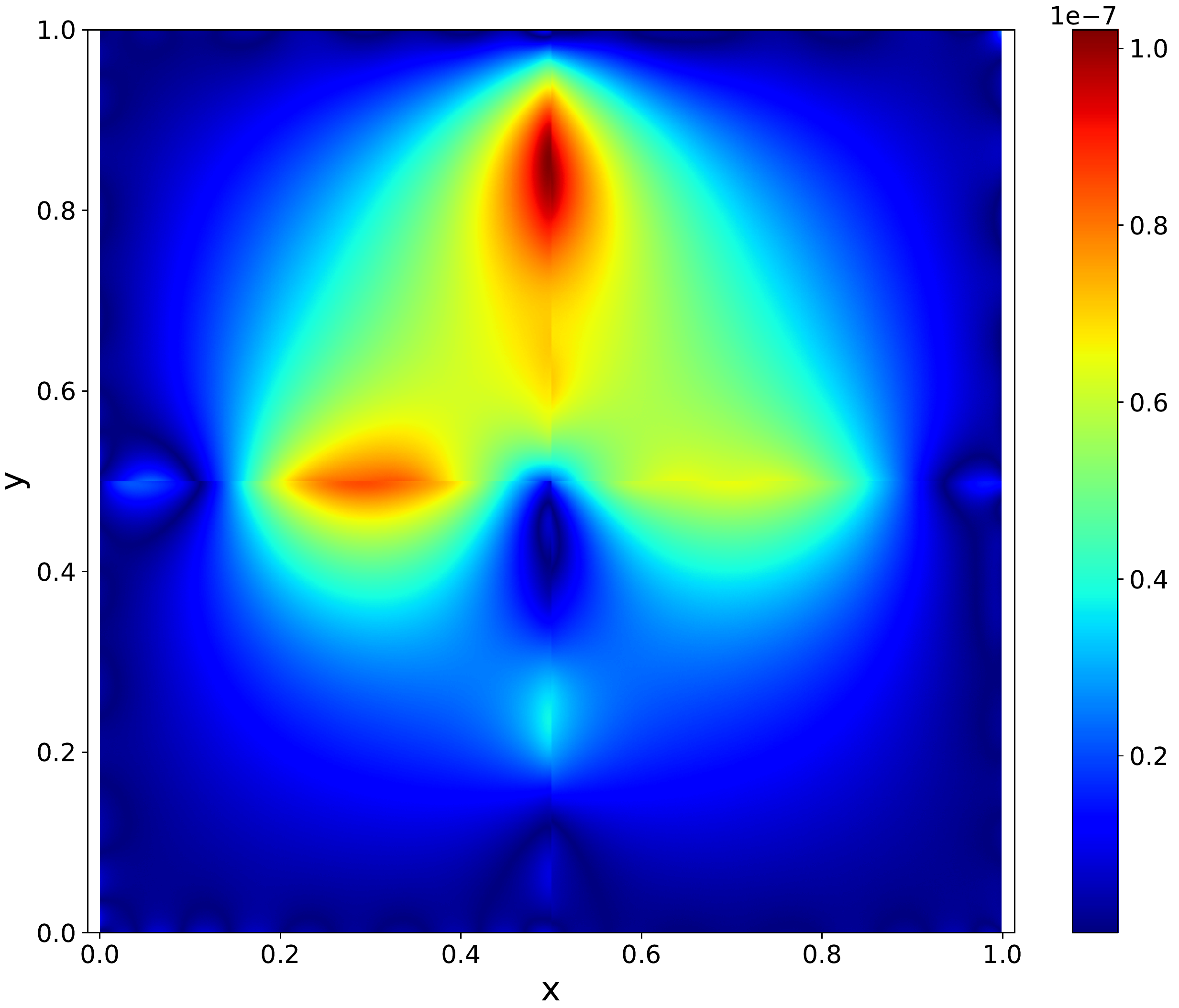}
			}
			\quad
			\subfigure[$\psi^{b}$]{
				\includegraphics[width=0.45\textwidth]{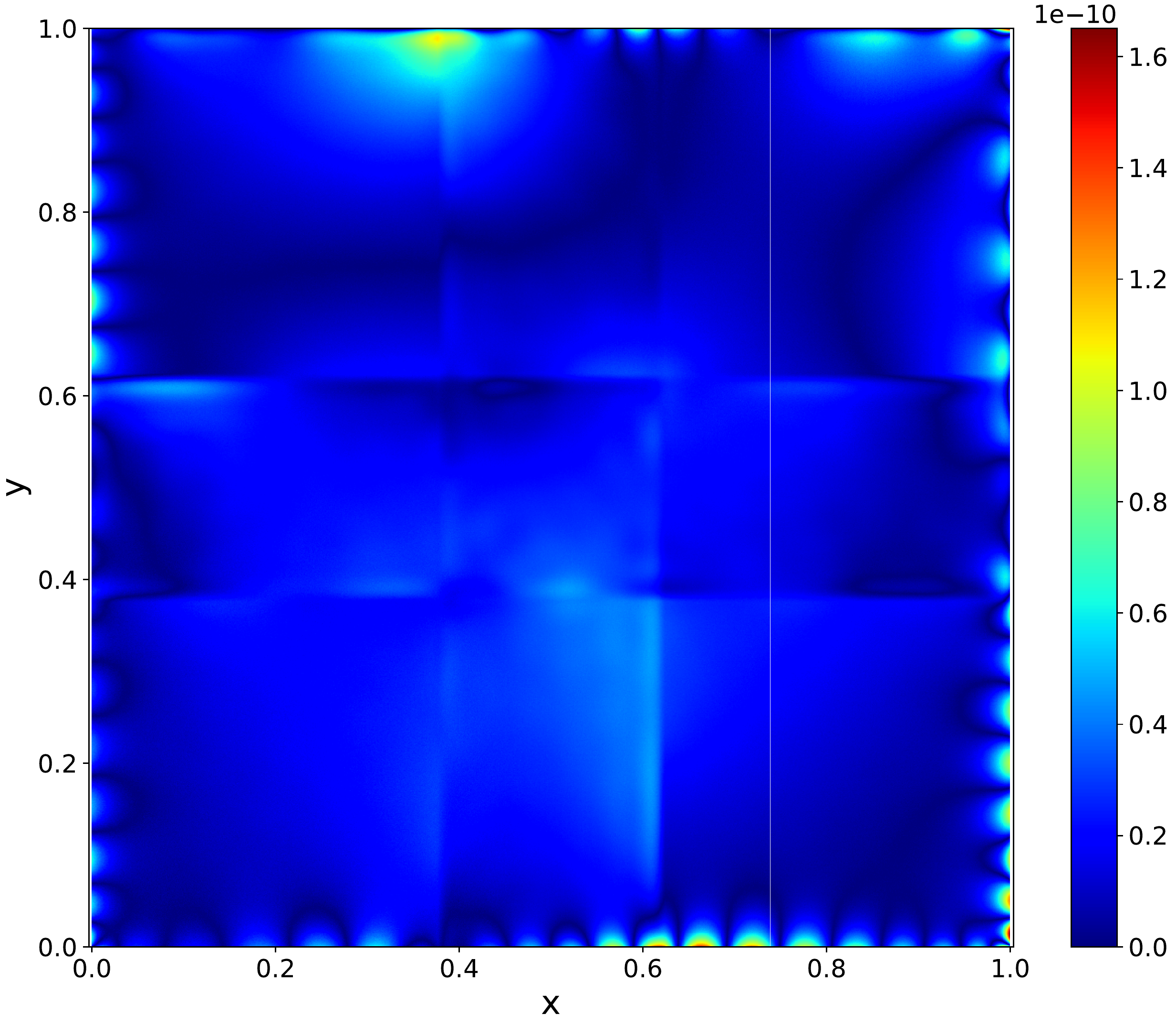}
			}
			\caption{Error distribution of the RFM with different choices of PoU for Poisson equation with solution \eqref{poisson-sol1}.}
			\label{fig3-1-1-2}
		\end{figure}

		Figure \ref{fig3-1-1-2} shows the error distribution of RFM for different choices of the PoU functions when $M=1600$ and $Q=2500$. For $\psi^{a}$, the error is more concentrated near the intersection of different sub-domains where smoothness conditions are imposed. For $\psi^{b}$, however, the error is concentrated near the boundary. Similar results are observed for the one-dimensional Helmholtz equation.
		
		\subsection{Multi-scale basis}
		\label{sec3-1-2}
		In this subsection, we show that the combination of local and global random feature functions works better in cases
		when the solution has both significant low and high frequency components.
		
		Consider the Poisson equation with the following explicit solution
		\begin{multline}
		u(x, y)=-A[\frac{3}{2} \cos (\pi x+\frac{2\pi}{5})+2 \cos (2 \pi x-\frac{\pi}{5})][\frac{3}{2} \cos (\pi y+\frac{2\pi}{5})+2 \cos (2 \pi y-\frac{\pi}{5})]\\
		-B[\frac{3}{2} \cos (2\pi x+\frac{4\pi}{5})+2 \cos (4 \pi x-\frac{2\pi}{5})][\frac{3}{2} (\cos 2\pi y+\frac{4\pi}{5})+2 \cos (4 \pi y-\frac{2\pi}{5})].
		\label{poisson-sol2}
		\end{multline}
		Three cases are considered: 1. a low-frequency problem when $A=1.0$, $B=0.0$; 2. a high-frequency problem when $A=0.0$, $B=1.0$; 3. mixed-frequency problem when $A=0.5$, $B=0.5$.
		
		The other hyper-parameters are set as follows:
		\begin{itemize}
			\item $M=1200$, $2700$, $4800$;
			
			\item $J_{n}=300$;
			
			\item $M_p=\frac{M}{300}$;
			
			\item $\boldsymbol{x}_{n} = (x_i, y_j)=(\frac{2i-1}{2\sqrt{M_p}},\frac{2j-1}{2\sqrt{M_p}}), i,j=1,\cdots,\sqrt{M_p}$;
			
			\item $\boldsymbol{r}_{n} = (r_i, r_j)=(\frac{1}{2\sqrt{M_p}},\frac{1}{2\sqrt{M_p}}), i,j=1,\cdots,\sqrt{M_p}$;
			
			\item $Q=1600$, $3600$, $6400$;
			
			\item $N=1920$, $4320$, $7680$.
		\end{itemize}
		Multi-scale basis functions are constructed using $\frac{300M_p}{M_p+1}$ basis functions for 
		the  values of $\{(x_i, y_j)\}$ and $\{(r_i, r_j)\}$ given above, together with $\frac{300M_p}{M_p+1}$ global basis functions with the parameters
		$(x, y) = (\frac{1}{2},\frac{1}{2})$, $(r_x, r_y)=(\frac{1}{2},\frac{1}{2})$. 
		$300$ local basis functions are constructed for each $(x_i, y_j)$ with the values of $(r_i, r_j)$ given above.
		
		\begin{table}[htbp]
			\caption{\label{table3-1-2-2} Comparison of PoU-based local basis and multi-scale basis functions for Poisson equation with the explicit solution \eqref{poisson-sol2}.} \centering
			\begin{small}
				\begin{tabular}{|c|c|c|c|c|}
					\hline
					Solution frequency & M & N & PoU-based basis & Multi-scale basis \\
					
					\hline
					\multirow{3}{*}{Low} &
					1200 & 1920 & 1.93E-8  & 3.28E-9  \\
					& 2700 & 4320 & 3.62E-9  & 6.42E-10 \\
					& 4800 & 7680 & 8.61E-10 & 3.05E-10 \\
					
					\hline
					\multirow{3}{*}{High} &
					1200 & 1920 & 6.42E-6 & 9.36E-7 \\
					& 2700 & 4320 & 1.34E-7 & 3.58E-8 \\
					& 4800 & 7680 & 4.16E-8 & 1.75E-8 \\
					
					\hline
					\multirow{3}{*}{Mixed} &
					1200 & 1920 & 3.22E-6 & 4.68E-7 \\
					& 2700 & 4320 & 6.54E-8 & 1.80E-8 \\
					& 4800 & 7680 & 2.06E-8 & 8.92E-9 \\
					
					\hline
				\end{tabular}
			\end{small}
		\end{table}
		
		Table \ref{table3-1-2-2} shows the error for RFM with multi-scale basis and with only PoU-based local basis functions. It is clear that the inclusion of global basis functions improves the accuracy when the solution has a significant low-frequency component. 
		A Fourier analysis of the errors confirms that the inclusion of global basis functions reduces the  low-frequency error more effectively.
		
		\subsection{Adaptive basis}
		\label{sec3-1-3}
		Here we demonstrate how the \textit{a prior} information helps us to select better random feature functions.
		
		Consider the one-dimensional Helmholtz equation with the following explicit solution
		\begin{equation}
		\begin{aligned}
		u(x) = & 4 \cos (4 (x+\frac{3}{20}))+5 \sin (\sqrt{5} (x+\frac{7}{20}))\\
		&+2 \sin (\sqrt{3} (x+\frac{1}{20}))+3 \sin(x+\frac{17}{20})+2.
		\end{aligned}
		\label{helmholtz-sol2}
		\end{equation}
		
		The other hyper-parameters are as follows:
		\begin{itemize}
			\item $M=400$, $800$, $1600$;
			
			\item $J_{n}=100$;
			
			\item $M_p=\frac{M}{100}$;
			
			\item $\boldsymbol{x}_{n} = x_n= 8\frac{2n-1}{2M_p}, n=1,\cdots,M_p$;
			
			\item $\boldsymbol{r}_{n} = r_{n}=\frac{8}{2M_p}, n=1,\cdots,M_p$;
			
			\item $Q=200$, $400$, $800$;
			
			\item $N=208$, $416$, $832$.
		\end{itemize}
		We use $\tanh$ and $\sin$ as the activation function. We test two different initialization methods: random initialization with 
		the distribution $\mathbb{U}[-R_{m},R_{m}]$ and equally spaced grids over $[-R_{m},R_{m}]$ with $\{k_m\}$ and $\{b_m\}$ being $-R_{m}+ 2R_{m}\frac{i}{10}, i=1,\cdots, 10$. Eight values of $R_{m}$. By the spectral analysis of the source term, we find that the highest frequency is $4$. Since a normalization is applied from $x$ to $\tilde{x}$, the highest-frequency term in $u$ corresponds to $k=\frac{4}{1/r_{m}}=\frac{16}{M_{p}}$. Therefore, the best results are obtained when the $\sin$ activation function is used over $[-R_{m},R_{m}]$ with $R_m\geq k$.
		
		Results of using adaptive random feature functions for this problem
		are shown in Table \ref{table3-1-3-1}. Another interesting observation is that 
		while choosing equally spaced feature vectors works for some cases, 
		the random initialization is found to be generally more reliable.
		\begin{table}[htbp]
			\caption{\label{table3-1-3-1} Results of the adaptive RFM for one-dimensional Helmholtz equation with solution \eqref{helmholtz-sol2}.} \centering
			\begin{small}
				\begin{tabular}{|c|c|c|c|c|c|}
					\hline
					\multirow{2}{*}{$M$} & \multirow{2}{*}{$R_{m}$} & \multicolumn{2}{|c|}{$\tanh$} & \multicolumn{2}{|c|}{$\sin$} \\
					\cline{3-6}
					& & $\mathbb{U}[-R_{m},R_{m}]$ & Equally spaced & $\mathbb{U}[-R_{m},R_{m}]$ & Equally spaced \\
					\hline
					\multirow{8}{*}{400}
					&1 & 3.03E-10 & 1.01E-10 & 6.21E-2  & 2.70E-2 \\
					&2 & 3.41E-11 & 1.75E-10 & 9.34E-5  & 1.23E-2 \\
					&3 & 5.12E-9  & 6.55E-10 & 4.68E-8  & 2.29E-3 \\
					&4 & 1.45E-7  & 2.93E-7  & 7.55E-13 & 8.02E-6 \\
					&5 & 1.77E-5  & 5.13E-4  & 1.14E-13 & 1.67E-4 \\
					&6 & 1.44E-4  & 3.02E-3  & 1.64E-13 & 9.60E-3 \\
					&7 & 8.20E-4  & 1.54E-2  & 2.31E-13 & 7.21E-2 \\
					&8 & 1.97E-2  & 8.58E-1  & 7.02E-14 & 5.62E-1 \\
					
					\hline
					\multirow{8}{*}{800}
					&1 & 9.78E-11 & 1.25E-11 & 2.44E-7  & 9.01E-6 \\
					&2 & 2.35E-11 & 1.66E-11 & 4.39E-13 & 3.77E-8 \\
					&3 & 2.11E-9  & 5.13E-11 & 2.52E-13 & 3.46E-6 \\
					&4 & 1.16E-7  & 1.01E-8  & 8.92E-13 & 1.04E-5 \\
					&5 & 3.03E-6  & 2.71E-5  & 1.02E-12 & 2.45E-4 \\
					&6 & 8.75E-5  & 2.22E-4  & 1.60E-12 & 3.46E-3 \\
					&7 & 5.94E-4  & 1.57E-3  & 2.17E-13 & 6.17E-2 \\
					&8 & 1.91E-3  & 9.51E-2  & 1.36E-12 & 5.12E-1 \\
					
					\hline
					\multirow{8}{*}{1600}
					&1 & 5.76E-12 & 7.29E-13 & 1.12E-12 & 1.68E-11 \\
					&2 & 6.64E-12 & 1.26E-11 & 6.23E-13 & 8.88E-7  \\
					&3 & 1.82E-9  & 1.50E-10 & 2.09E-13 & 3.48E-5  \\
					&4 & 2.63E-7  & 7.59E-9  & 2.04E-13 & 2.18E-4  \\
					&5 & 6.00E-6  & 6.25E-6  & 1.29E-12 & 3.98E-3  \\
					&6 & 1.22E-4  & 8.59E-5  & 3.96E-12 & 1.07E-2  \\
					&7 & 4.55E-3  & 7.24E-4  & 1.16E-12 & 2.39E-1  \\
					&8 & 4.65E-3  & 5.16E-2  & 9.79E-13 & 2.07E+0  \\
					\hline
				\end{tabular}
			\end{small}
		\end{table}
		
		We now turn to the Poisson equation with the solution \eqref{poisson-sol1}.
		We set the hyper-parameters as follows:
		\begin{itemize}
			\item $M=4000$;
			
			\item $J_{n}=1000$;
			
			\item $M_p=\frac{M}{1000}$;
			
			\item $\boldsymbol{x}_{n} \in \left\{(\frac{1}{4},\frac{1}{4}), (\frac{1}{4},\frac{3}{4}), (\frac{3}{4},\frac{1}{4}), (\frac{3}{4},\frac{3}{4})\right\}$;
			
			\item $\boldsymbol{r}_{n} =(r_x, r_y)=(\frac{1}{4},\frac{1}{4})$;
			
			\item $Q=1600$;
			
			\item $N=1920$.
		\end{itemize}
		We use $\tanh$ and $\sin$ as the activation function. We test two different initialization methods: random initialization with distribution $\mathbb{U}[-R_{m},R_{m}]$ and equally spaced grids over $[-R_{m},R_{m}]$ with $\{k^1_m, k^2_m\}$ and $\{b_m\}$ being $-R_{m}+ 2R_{m}\frac{i}{10}, i=1,\cdots, 10$. $10$ values of $R_{m}$ are tested. Results of using adaptive random feature functions 
		are shown in Table \ref{table3-1-3-2}. Again, the best results are observed when the $\sin$ activation function is used over $[-R_{m},R_{m}]$ with $R_m\geq k$ and random initialization is found to be generally more reliable.	
		
		\begin{table}[htbp]
			\caption{\label{table3-1-3-2} Results of using adaptive random feature functions for the two-dimensional Poisson equation with solution \eqref{poisson-sol1}.} 
			\centering
			\begin{small}
				\begin{tabular}{|c|c|c|c|c|}
					\hline
					\multirow{3}{*}{$R_{m}$} & \multicolumn{2}{|c|}{$\tanh$} & \multicolumn{2}{|c|}{$\sin$} \\
					\cline{2-5}
					& $\mathbb{U}[-R_{m},R_{m}]$ & Equally spaced & $\mathbb{U}[-R_{m},R_{m}]$ & Equally spaced\\
					\hline
					0.5 & 4.92E-9   & 1.01E-9   & 2.55E-3   & 6.05E-4  \\
					1.0 & 2.91E-8   & 9.36E-9   & 8.96E-7   & 2.58E-5  \\
					1.5 & 1.33E-6   & 5.95E-7   & 1.79E-9   & 1.47E-6  \\
					2.0 & 8.75E-5   & 7.85E-5   & 3.30E-12  & 4.29E-7  \\
					2.5 & 8.16E-4   & 4.70E-5   & 2.86E-12  & 7.66E-6  \\
					3.0 & 2.06E-2   & 5.27E-4   & 7.32E-12  & 2.17E-5  \\
					3.5 & 1.53E-3   & 3.95E-3   & 6.10E-12  & 7.45E-5  \\
					4.0 & 2.66E-3   & 1.27E-3   & 6.10E-12  & 5.59E-5  \\
					4.5 & 5.39E-3   & 1.76E-2   & 2.29E-11  & 1.24E-3  \\
					5.0 & 1.29E-2   & 5.16E-2   & 2.17E-11  & 6.72E-3  \\
					\hline
				\end{tabular}
			\end{small}
		\end{table}

		\section{Experimental setup and numerical results for  problems with complex geometry}
		\label{sec::appendix::complex}
		\subsection{Rescaling}
		\label{sec::appendix::complex::rescaling}
		The exact displacement solution for the Timoshenko beam problem is
		\begin{equation}
		\begin{aligned}
		&u=-\frac{P y}{6 E I}[(6 L-3 x) x+(2+\nu)(y^{2}-\frac{D^{2}}{4})], \\
		&v=\frac{P}{6 E I}[3 \nu y^{2}(L-x)+(4+5 \nu) \frac{D^{2} x}{4}+(3 L-x) x^{2}],
		\end{aligned}
		\label{timoshenko-sol}
		\end{equation}
		where $I = \frac{D^{3}}{12}$. Homogeneous Dirichlet boundary condition is applied on the left boundary $x=0$ and 
		Homogeneous Neumann boundary condition is applied on the other boundaries. The material parameters are as follows: the Young's modulus $E=3 \times 10^{7}$ Pa, Poisson ratio $\nu=0.3$. We choose $D = 10$, $L = 10$, and the shear force is $P = 1000$ Pa. 
		
		The other hyper-parameters for the Timoshenko beam problem are as follows:
		\begin{itemize}
			\item $M=800$, $3200$;
			
			\item $J_{n}=200$;
			
			\item $M_p=\frac{M}{200}$;
			
			\item $\boldsymbol{x}_{n} = (x_i, y_j)=(10\frac{2i-1}{2\sqrt{M_p}},10\frac{2j-1}{2\sqrt{M_p}}), i,j=1,\cdots,\sqrt{M_p}$;
			
			\item $\boldsymbol{r}_{n} = (r_x, r_y)=(5,5)$;
			
			\item $Q=25M_p,100M_p,400M_p,1600M_p$.
		\end{itemize}
		We construct $\frac{200M_p}{M_p+1}$ basis functions associated with the choice of  $\{(x_i, y_j) \}$ with $\{(r_x, r_y)\}$
		given above, and add $\frac{200M_p}{M_p+1}$ basis functions associated with the point $(5,5)$ with $(r_x, r_y)=(5,5)$. To count the total number of basis functions and the number of conditions, we convert $N_{x}$, $N_{y}$, $Q_{x}$, $Q_{y}$, and $M^{\prime}$ in locELM to $M$ and $N$ in RFM according to $M = 2 N_{x}N_{y} M^{\prime}$ and
		$
		N = 2 N_{x}N_{y}Q_{x}Q_{y} + 4N_{x}Q_{x} + 4N_{y}Q_{y} + 6N_{x}Q_{x}(N_{y}-1) + 6N_{y}Q_{y}(N_{x}-1).
		$
		
		The exact displacement field for the two-dimensional elasticity problem with complex geometry shown in Figure \ref{fig3-2-2} is
		\begin{equation}
		\begin{aligned}
		&u=\frac{1}{10} y ((x+10)  \sin y+(y+5)  \cos x), \\
		&v=\frac{1}{60} y ((30+5 x \sin (5 x)) (4+e^{-5 y})-100).
		\end{aligned}
		\label{2d-complex-sol}
		\end{equation}
		Dirichlet boundary condition is applied on the lower boundary $y=0$ and Neumann boundary condition is applied on the other boundaries and the holes inside.
		The material constants are:  the Young's modulus $E=3 \times 10^{7}$ Pa and Poisson ratio $\nu=0.3$. 
		
		The other hyper-parameters are as follows:
		\begin{itemize}
			\item $M=3200$, $12800$;
			
			\item $J_{n}=200$;
			
			\item $M_p=\frac{M}{200}$;
			
			\item $\boldsymbol{x}_{n} = (x_i, y_j)=(8\frac{2i-1}{2\sqrt{M_p}},8\frac{2j-1}{2\sqrt{M_p}}), i,j=1,\cdots,\sqrt{M_p}$;
			
			\item $\boldsymbol{r}_{n} = (r_x, r_y)=(\frac{8}{2\sqrt{M_p}},\frac{8}{2\sqrt{M_p}})$;
			
			\item $Q=25M_p$, $100M_p$, $400M_p$, $1600M_p$.
		\end{itemize}
		We construct $\frac{200M_p}{M_p+1}$ basis functions for each point $(x_i, y_j)$ with $(r_x, r_y)$ given above, and adds $\frac{200M_p}{M_p+1}$ basis functions for the point $(4,4)$ with $(r_x, r_y)=(4,4)$.
		
		Results of the RFM for the elasticity problem with this explicit solution are shown in Table \ref{table3-2-2}.
		\begin{table}[htbp]
			\caption{\label{table3-2-2} Results of RFM for the elasticity problem with  complex geometry.} \centering
			\begin{small}
				\begin{tabular}{|c|cc|ccccc|}
					\hline
					Method  & $M$ & $N$ & $u$ error  & $v$ error & $\sigma_{x}$ error & $\sigma_{y}$ error & $\tau_{xy}$ error \\
					
					\hline
					\multirow{8}{*}{RFM}
					&  \multirow{4}{*}{3200} & 1784   & 4.96E-1  & 8.37E-1  & 1.09E+0  & 3.52E+0 & 5.24E-1  \\
					&                        & 4658   & 5.82E-3  & 7.12E-3  & 1.04E-2  & 5.47E-2 & 3.85E-3  \\
					&                        & 13338  & 1.69E-5  & 1.19E-5  & 2.89E-5  & 6.40E-5 & 8.18E-6  \\
					&                        & 42820  & 1.39E-5  & 1.55E-5  & 4.92E-5  & 6.16E-5 & 1.29E-5  \\
					\cline{2-8}
					&  \multirow{4}{*}{12800}& 6578   & 9.11E-2  & 6.41E-2  & 1.03E-1  & 2.46E-1 & 2.95E-2  \\
					&                        & 17178  & 2.35E-4  & 2.10E-4  & 3.02E-4  & 7.56E-4 & 8.93E-5  \\
					&                        & 50500  & 5.46E-7  & 4.98E-7  & 8.45E-7  & 2.03E-6 & 2.67E-7  \\
					&                        & 165184 & 2.32E-7  & 1.89E-7  & 9.28E-8  & 2.32E-7 & 2.43E-8  \\
					\hline
				\end{tabular}
			\end{small}
		\end{table}
		
		\subsection{Comparison with FEM}
		\label{sec::appendix::complex::fem}
		For the two-dimensional elasticity problem in Section \ref{sec3-3}, we set hyper-parameters as follows:
		\begin{itemize}
			\item $M=16000$;
			
			\item $J_{n}=400$;
			
			\item $M_p=\frac{M}{400}$;
			
			\item $\boldsymbol{x}_{n} = (x_i, y_j) = (\frac{2i-1}{8}-1,\frac{2j-1}{8}-\frac{1}{2}), i=1,\cdots,10, j=1,\cdots,4$;
			
			\item $\boldsymbol{r}_{n} = (r_x, r_y)=(\frac{1}{8},\frac{1}{8})$;
			
			\item $Q=16000$, $64000$, $144000$, $256000$;
			
			\item $N=40326$, $135442$, $285472$, $490176$.
		\end{itemize}
		For comparison, we implement the standard adaptive FEM with total degrees of freedom $M=3716$, $10438$, $40054$, $153562$.
		
		Table \ref{table3-3-1} shows the error between RFM and FEM for the elasticity problem in Section \ref{sec3-3}.
		\begin{table}[htbp]
			\caption{\label{table3-3-1} Comparison of the numerical solutions of RFM and FEM for the elasticity problem.} \centering
			\begin{small}
				\begin{tabular}{|c|c|cc|ccccc|}
					\hline
					Method & Reference & $M$ & $N$ & $u$ error  & $v$ error & $\sigma_{x}$ error & $\sigma_{y}$ error & $\tau_{xy}$ error \\
					
					\hline
					\multirow{3}{*}{RFM} & 			\multirow{3}{*}{RFM $N=490176$}
					& \multirow{3}{*}{16000}& 40326  &  1.28E+0  & 1.12E+0  & 1.29E+0  & 9.37E-1 & 1.03E+0  \\
					& &  & 135442 &  1.12E-1  & 1.16E-1  & 1.13E-1  & 1.03E-2 & 1.20E-1  \\
					& &  & 285472 &  6.52E-4  & 6.98E-4  & 1.03E-3  & 3.01E-5 & 1.88E-3  \\
					
					\hline
					\multirow{4}{*}{RFM} & \multirow{4}{*}{FEM $M=153562$}
					& \multirow{4}{*}{16000}  & 40326  &  1.30E+0  & 1.12E+0  & 1.28E+0  & 9.37E-1 & 1.03E+0  \\
					& &  & 135442 &  7.65E-2  & 8.55E-2  & 1.16E-1  & 1.31E-1 & 1.25E-1  \\
					& &  & 285472 &  3.94E-2  & 3.36E-2  & 6.59E-3  & 5.95E-2 & 2.31E-2  \\
					& &  & 490176 &  4.00E-2  & 3.43E-2  & 6.20E-3  & 5.92E-2 & 2.30E-2  \\
					
					\hline
					\multirow{3}{*}{FEM} & \multirow{3}{*}{FEM $M=153562$}
					& 3716    & 3716 &  3.15E-4  & 4.54E-4  & 1.41E-2  & 5.81E-2 & 3.35E-2  \\
					& & 10438   & 10438 &  1.20E-4  & 1.81E-4  & 9.39E-3  & 3.61E-2 & 2.13E-2  \\
					& & 40054   & 40054 &  2.88E-5  & 3.93E-5  & 4.65E-3  & 1.62E-2 & 9.40E-3  \\
					
					\hline
					\multirow{4}{*}{FEM} & \multirow{4}{*}{RFM $N=490176$}
					& 3716    & 3716 &  3.87E-2  & 3.36E-2  & 1.43E-2  & 8.93E-2 & 3.86E-2  \\
					& & 10438   & 10438 &  3.86E-2  & 3.34E-2  & 1.05E-2  & 7.29E-2 & 2.99E-2  \\
					& & 40054   & 40054 &  3.85E-2  & 3.32E-2  & 7.19E-3  & 6.33E-2 & 2.44E-2  \\
					& & 153562  & 153562 &  3.85E-2  & 3.32E-2  & 6.22E-3  & 6.01E-2 & 2.31E-2  \\
					\hline
				\end{tabular}
			\end{small}
		\end{table}
		
		We set the hyper-parameters for the elasticity problem over a complex geometry in Section \ref{sec3-3} as follows:
		\begin{itemize}
			\item $M=14400$;
			
			\item $J_{n}=400$;
			
			\item $M_p=\frac{M}{400}$;
			
			\item $\boldsymbol{x}_{n} = (x_i, y_j) = (8\frac{2i-1}{12},8\frac{2j-1}{12}), i,j=1,\cdots,6$;
			
			\item $\boldsymbol{r}_{n} = (r_x, r_y)=(\frac{8}{12},\frac{8}{12})$;
			
			\item $Q=129600$, $152100$, $176400$, $202500$, $230400$;
			
			\item $N=195146$, $226132$, $259400$, $294878$, $332606$.
		\end{itemize}
		
		Table \ref{table3-3-2} records the results of RFM for the  elasticity problem over a complex geometry in Section \ref{sec3-3}.
		\begin{table}[htbp]
			\caption{\label{table3-3-2} Numerical results of the RFM for the elasticity problem over a complex geometry. The result with  $N=332606$ is taken as the reference solution.} \centering
			\begin{small}
				\begin{tabular}{|c|c|ccccc|}
					\hline
					$M$ & $N$ & $u$ error  & $v$ error & $\sigma_{x}$ error & $\sigma_{y}$ error & $\tau_{xy}$ error \\
					
					\hline
					\multirow{4}{*}{$14400$}
					& 195146  &  2.30E-1  & 1.30E-1  & 6.64E-2  & 1.72E-1 & 1.71E-1  \\
					& 226132  &  8.97E-2  & 1.23E-1  & 5.60E-2  & 1.41E-1 & 1.32E-1  \\
					& 259400  &  6.47E-2  & 6.94E-2  & 3.66E-2  & 9.04E-2 & 8.15E-2  \\
					& 294878  &  7.30E-2  & 6.68E-2  & 3.46E-2  & 7.13E-2 & 7.05E-2  \\
					\hline
				\end{tabular}
			\end{small}
		\end{table}

		\subsection{The elliptic homogenization problem}
		\label{sec::appendix::complex::homogenization}
		For the homogenization problem, we set the hyper-parameters as follows:
		\begin{itemize}
			\item $M=25600$;
			
			\item $J_{n}=400$;
			
			\item $M_p=\frac{M}{400}$;
			
			\item $\boldsymbol{x}_{n} = (x_i, y_j) =  (\frac{2i-1}{8}-1,\frac{2j-1}{8}-1), i, j = 1,\cdots,8$;
			
			\item $\boldsymbol{r}_{n} = (r_x, r_y)=(\frac{1}{8},\frac{1}{8})$;
			
			\item $Q=25600$, $102400$, $230400$, $409600$;
			
			\item $N=25554$, $91339$, $197360$, $343586$.
		\end{itemize}
		
		Results of RFM for the homogenization problem in Section \ref{sec3-4} are shown in Table \ref{table3-4}.
		\begin{table}[htbp]
			\caption{\label{table3-4} Numerical convergence of the  random feature method for the homogenization problem.} \centering
			\begin{small}
				\begin{tabular}{|c|c|ccc|}
					\hline
					$M$ & $N$ & $u$ error & $u_x$ error & $u_y$ error \\
					\hline
					\multirow{4}{*}{25600}
					& 25554    &  1.42E+0   &8.68E+0  &8.73E+0 \\
					& 91339    &  3.13E-2   &3.54E-2  &3.62E-2 \\
					& 197360   &  3.48E-3   &6.45E-3  &7.18E-3 \\
					& 343586   &  \multicolumn{3}{|c|}{Reference}  \\
					\hline
				\end{tabular}
			\end{small}
		\end{table}

		\subsection{Stokes flow}
		\label{sec::appendix::complex::stokes}
		The exact displacement field for the Stokes flow is given by
		\begin{equation}
		\begin{aligned}
		u=&x+x^{2}-2 x y+x^{3}-3 x y^{2}+x^{2} y, \\
		v=&-y-2 x y+y^{2}-3 x^{2} y+y^{3}-x y^{2},\\
		p=&x y+x+y+x^{3} y^{2}-\frac{4}{3}.
		\end{aligned}
		\end{equation}
		
		We set the hyper-parameters for the Stokes flow as follows:
		\begin{itemize}
			\item $M=400$, $800$, $1600$;
			
			\item $J_{n}=100$, $200$, $400$;
			
			\item $M_p= 4$;
			
			\item $\boldsymbol{x}_{n} \in \left\{(\frac{1}{4},\frac{1}{4}),(\frac{1}{4},\frac{3}{4}),(\frac{3}{4},\frac{1}{4}),(\frac{3}{4},\frac{3}{4})\right\}$;
			
			\item $\boldsymbol{r}_{n} = (r_x, r_y)=(\frac{1}{4},\frac{1}{4})$;
			
			\item $Q=100$, $400$, $1600$, $6400$;
			
			\item $N=512$, $1596$, $5390$, $19488$.
		\end{itemize}

		Results of the RFM are shown in Table \ref{table3-5-1}.
		\begin{table}[htbp]
			\caption{\label{table3-5-1} Numerical results of the RFM for the Stokes flow with an explicit solution.} \centering
			\begin{small}
				\begin{tabular}{|c|c|ccc|}
					\hline
					$M$ & $N$ & $u$ error & $v$ error & $p$ error\\
					\hline
					\multirow{4}{*}{400}
					& 512   & 3.22E-4 & 2.28E-4 & 3.21E-2 \\
					& 1596  & 6.13E-7 & 3.44E-7 & 9.72E-5 \\
					& 5390  & 4.22E-7 & 2.54E-7 & 1.64E-4 \\
					& 19488 & 1.44E-7 & 1.03E-7 & 1.31E-5 \\
					\hline
					\multirow{4}{*}{800}
					& 512   & 5.25E-4  & 3.49E-4  & 4.39E-2 \\
					& 1596  & 4.95E-7  & 3.03E-7  & 2.77E-5 \\
					& 5390  & 1.60E-10 & 9.48E-11 & 1.73E-7 \\
					& 19488 & 1.15E-10 & 6.01E-11 & 1.06E-7 \\
					\hline
					\multirow{4}{*}{1600}
					& 512   & 3.33E-4  & 3.61E-4  & 4.64E-2 \\
					& 1596  & 1.11E-6  & 6.23E-7  & 5.67E-5 \\
					& 5390  & 3.02E-12 & 1.56E-12 & 1.06E-9 \\
					& 19488 & 2.45E-13 & 1.63E-13 & 1.37E-9 \\
					\hline
				\end{tabular}
			\end{small}
		\end{table}
		
	\end{appendix}
	
	\bibliographystyle{siam}
	\bibliography{RFM}
	
\end{document}